\documentclass[letterpaper,12pt]{amsart}

\usepackage{amscd,amssymb,pb-diagram,fullpage}

\newcommand{\cale}{{\mathcal E}}

\newcommand{\p}{{\mathfrak p}}
\newcommand{\q}{{\mathfrak q}}

\newcommand{\oo}{{\mathcal O}}
\newcommand{\ok}{{{\mathcal O}_k}}
\newcommand{\C}{{\mathbf C}}
\newcommand{\F}{{\mathbf F}}
\newcommand{\fp}{{\mathbf F_p}}

\newcommand{\Q}{{\mathbf Q}}

\newcommand{\Z}{{\mathbf Z}}

\newcommand{\ef}{{\epsilon_f}}
\newcommand{\ea}{{\epsilon_a}}
\newcommand{\e}{{\epsilon}}

\newcommand{\vsp}{\vspace{8pt}}

\newcommand{\rarr}{{\rightarrow}}

\newcommand{\mymod}[2]{{ #1 \: (\bmod \: {#2})}}

\newtheorem{thm}{Theorem}
\newtheorem{lem}{Lemma}

\newtheorem{cor}{Corollary}

\theoremstyle{remark}
\newtheorem{quest}{Question}
\newtheorem{remark}{Remark}

\newtheorem*{notat}{Notations}

\newtheorem*{convent}{Convention}

\newcommand{\et}{{\text{\'et}}}
\newcommand\gal{\text{Gal}}

\newcommand{\binomial}[2]{{
  \Big( \begin{array}{c} {#1} \\ {#2} \end{array} \Bigr)
}}

\newcommand{\ov}[1]{{\overline{{#1}}}}
\newcommand{\efc}{{f: \mathcal E\rarr C}}

\newcommand{\E}{{\mathcal E}}
\newcommand{\etwo}{{\mathcal E^{(2)}}}
\newcommand{\etwostar}{{\mathcal E^{(2)}_\ast}}
\newcommand{\etwotilde}{{\widetilde{\E}^{(2)}}}

\newcommand{\ej}{{\mathcal E^{(j)}}}
\newcommand{\ejtilde}{{\widetilde{\E}^{(j)}}}

\newcommand{\stwo}{{S^{(2)}}}
\newcommand{\stwotilde}{{\widetilde{S}^{(2)}}}

\DeclareFontFamily{U}{wncyr}{}
\DeclareFontShape{U}{wncyr}{m}{n}{%
   <5> <6> <7> <8> <9> gen * wncyr
   <10> <10.95> <12> <14.4> <17.28> <20.74>  <24.88>wncyr10}{}
\DeclareFontShape{U}{wncyr}{bx}{n}{%
   <5> <6> <7> <8> <9> gen * wncyb
   <10> <10.95> <12> <14.4> <17.28> <20.74>  <24.88>wncyb10}{}
\DeclareSymbolFont{cyr}{U}{wncyr}{m}{n}
\SetSymbolFont{cyr}{bold}{U}{wncyr}{bx}{n}

\DeclareMathSymbol{\sha}{\mathalpha}{cyr}{"58}

\font\twcmbx=cmbx12
\newfam\swfam  
\textfont\swfam=\twcmbx

\newcount\timehh\newcount\timemm
\timehh=\time 
\divide\timehh by 60 \timemm=\time
\begin{document}

\title{On the N\'eron-Severi groups of fibered varieties}

\author{Siman Wong
  {
    \protect \protect\sc\today\ -- 
    \ifnum\timehh<10 0\fi\number\timehh\,:\,\ifnum\timemm<10 0\fi\number\timemm
    \protect \, \, \protect \bf DRAFT
  }
}

\begin{abstract}
We apply Tate's conjecture on algebraic cycles to study
the N\'eron-Severi groups of varieties fibered over a curve.
This is inspired by the work of Rosen and Silverman, who carry out such
an analysis to 
derive a formula for the rank of the group of sections of an elliptic
surface.
For a semistable fibered surface, under Tate's conjecture we derive a
formula for the rank of the group of sections of the associated Jacobian
fibration.
For fiber powers of a semistable elliptic fibration
$
\E\rarr C
$, under Tate's conjecture
we give a recursive formula for the rank of the N\'eron-Severi groups of these
fiber powers.   For fiber squares, we construct 
unconditionally a set of
independent elements in the N\'eron-Severi groups.   When 
$
\E\rarr C
$
is the universal elliptic curve over the modular curve
$
X_0(M)/\Q
$,
we apply the Selberg trace formula to verify
our recursive formula in the case of fiber squares.
This gives an analytic proof of Tate's conjecture for such fiber squares
over $\Q$, and it shows that the independent elements we constructed in fact
form a basis of the N\'eron-Severi groups.
This is the fiber square analog of the Shioda-Tate Theorem.
\end{abstract}

\address{Department of Mathematics, University of Massachusetts.
        Amherst, MA 01003-4515}

\email{siman@math.umass.edu}





\maketitle

\tableofcontents

\section{Introduction}

Let $K$ be a field which is finitely generated over its prime field, and let
$A$ be an Abelian variety defined over $K$. The Mordell-Weil Theorem
\cite[p.~138]{lang} shows that $A(K)$ is a finitely generated Abelian
group.  It is of great number-theoretic interest to determine the rank of
$A(K)$.  If $K$ is a global field, the Birch-Swinnerton-Dyer
conjecture (as generalized by Tate \cite{tate:bsd}) predicts that the
Hasse-Weil $L$-function $L(A/K, s)$ of $A/K$ has an analytic continuation to
the whole complex plane, and that it has a zero at $s=1$ of order equal to the
rank of $A(K)$.  This latter statement is equivalent to saying that the rank
of $A(K)$ is equal to the residue at $s=1$ of the logarithmic derivative of
$L(A/K, s)$.
%
%
Now, take $K=\Q$ and $A$ to be an elliptic curve $E$ over $\Q$.
A
formal\footnote{the
	Tauberian theorem is not applicable since the $a_p(E)$'s change
	sign.}
Tauberian argument then leads Nagao \cite{nagao} to observe that the
expression
\begin{equation}
\frac{-1}{x}\sum_{p\le x} a_p(E)\log p,
	\label{expression}
\end{equation}
where
$
a_p(E) := p+1- \#(\text{non-singular points of $E$ mod $p$})
$,
approximates the rank of $E(\Q)$ for suitable, large
values of $x$.
%
%
This is supported by computer data \cite{mestre}.
For an elliptic curve $\cale$ defined over the rational function field
$\Q(T)$, Nagao
\cite{nagao}
defines a similar expression
\begin{equation}
S(\cale, x)
:=
\frac{-1}{x}\sum_{p\le x}\frac{1}{p}\sum_{t\in\F_p}a_p(\cale_t)\log p,
	\label{nagao}
\end{equation}
where $\mathcal \cale_t$ denotes the curve obtained from $\cale$ by
specializing the variable $T$ to $t$. For several families of elliptic curves
$\cale$ over
$\Q(T)$, Nagao \cite{nagao} shows that $S(\cale, x)$ converges to the
Mordell-Weil rank
of $\cale$ over $\Q(T)$, and he conjectures that this is true in general:
\begin{equation}
\text{MW-rank } \cale(\Q(T))
\stackrel{?}{=}
\lim_{x\rarr\infty}
	\frac{-1}{x} \sum_{p\le x}\frac{1}{p}\sum_{t\in\F_p}
	a_p(\cale_t)\log p.
	\label{nagao_conj}
\end{equation}

Elliptic curves over $\Q(T)$ can be viewed as elliptic surfaces over
$
\mathbf P^1_\Q
$.
More precisely, let  $C$ be a smooth, geometrically connected projective curve
defined over a number field $k$; denote by $K_C=k(C)$ its function field.
Let $E$ be an elliptic curve over $K_C$.  By the theory
of minimal models \cite{bosch}, there exists a smooth projective surface
$\E$ together with a genus one fibration
$
f: \E\rarr C
$
with a section $\sigma$, all defined over $k$, such that the generic fiber of
$f$ is
$
E/K_C
$,
and that no fiber of $f$ contains a curve
$
\simeq \mathbf P^1
$
and with self-intersection number $-1$; such elliptic fibrations are called 
relatively minimal.  Furthermore, the correspondence
$
E/K_C \leftrightarrow (\E, f)
$
is bijective \cite{licht}.  Denote by
$
NS(\E/k)
$
the N\'eron-Severi group of the surface $\E/k$.   Fix an algebraic closure
$\ov{k}$ of $k$, and write  $G_k$ for the absolute Galois group
$
\gal(\ov{k}/k)
$.
Given a maximal ideal $\p$ of the ring of integers $\oo_k$ of $k$, denote by
$N_\p$
the absolute norm of $\p$.
\if 3\
{
Tate's conjecture \cite{tate} on algebraic cycles predicts that the rank of
$
NS(\mathcal E/k)
$
is equal to the order of \textit{pole} at $s=2$ of the $L$-function attached to
$
H^2_{\text{\'et}}(\mathcal E\otimes\ov{k}, \Q_l(1))
$.
}
\fi
For any smooth variety $V/k$ of dimension $>0$,  Tate's conjecture on
algebraic cycles \cite{tate} states that 
$
L_2(V/k, s)
$,
the $L$-function attached to
$
H^2_\et(V\otimes\ov{k}, \Q_l)
$,
has meromorphic continuation to $\C$, and that
\begin{eqnarray}
&&
NS(V/k)\otimes\Q_l \simeq H^2_\et(V\otimes\ov{k}, \Q_l(1))^{G_k},
	\label{conj1}
\\
&&
-\raisebox{-7pt}{$\stackrel{\displaystyle\text{ord}}{\scriptstyle s=2}$}\,
	L_2(V/k, s)
=
\dim H^2_\et(V\otimes\ov{k}, \Q_l(1))^{G_k}.
	\label{conj2}
\end{eqnarray}
\if 3\
{
\begin{equation}

\simeq
H^2_{\text{\'et}}(\mathcal E\otimes\ov{k}, \Q_l(1))^{G_k}
	\label{tate}
\end{equation}
and
\begin{equation}
-\raisebox{-7pt}{$\stackrel{\displaystyle\text{ord}}{\scriptstyle s=2}$}\,
	L_2(\mathcal E/k, s)
=
\text{rank} NS(\mathcal E/k).
	\label{tate2}
\end{equation}
}
\fi
On the other hand, the Shioda-Tate formula \cite{shioda:mw}  furnishes a
natural isomorphism
\begin{equation}
	E(\ov{k}(C))\simeq NS(\cale/{\ov{k}})/T_\cale,
\label{shioda-tate}
\end{equation}
where $T_\cale$ denotes the subgroup of $NS(\cale/{\ov{k}})$ generated by the
section $\sigma$ together with all the $\ov{k}$-components of the fibers of
$f$.  This leads Silverman to consider the following analytic form of Nagao's
conjecture
\begin{eqnarray}
\text{MW-rank of } \cale( k(C) )
\stackrel{?}{=}
\raisebox{-7pt}{$\stackrel{\displaystyle\text{res}}{\scriptstyle s=1}$}
\,
	\sum_{\p}\frac{-\log N_\p}{N_{\p}^{s+1}}\sum_{t\in C(\F_\p)}
	a_{\p}(\cale_t)
\label{analy}
\end{eqnarray}
and asks if it
is related to Tate's conjecture. By giving an
$L$-series interpretation of the Shioda-Tate formula and by analyzing
carefully the singular fibers of the elliptic surface $\cale$, Rosen and
Silverman
\cite{mike_joe} show that, for any non-split elliptic surface
$
\efc
$
with a section,
(\ref{analy}) indeed follows from Tate's conjecture.
%
%
%
%
The original Nagao's conjecture
then follows from (\ref{analy}) plus a non-vanishing hypothesis for
$L_2(\mathcal E/\Q, s)$.  In particular, Rosen and Silverman show that
Nagao's conjecture holds unconditionally for rational elliptic surfaces
with sections.

This result of Rosen and Silverman can be thought of as a
Birch-Swinnerton-Dyer-type  conjecture for \textit{families} of elliptic
curves over $C/k$.  If $\E_t$ has good reduction at $\p$, then 
$
a_\p(\E_t)
$
is the trace of the geometric Frobenius at $\p$ of
$
H^1_\et(C\otimes\ov{k}, \Q_l)
$.
A natural question then arises:
are there analogs of Nagao's conjecture for general fibered varieties over
a curve?   In this paper we investigate this question for fibered surfaces
and for fiber powers of elliptic fibrations.

\if 3\
{
In this paper we consider the analytic form of Nagao's conjecture for two
types of fibered varieties.  For a wide class of
fibered surface (including the semistable ones) we show that under the Tate
conjecture, the analog of (\ref{analy}) again gives the rank of the group of
	sections.   
For fiber powers of a semistable elliptic fibration
$
\E\rarr C
$, under the Tate conjecture
we give a recursive formula for the rank of the N\'eron-Severi groups of these
fiber powers.   In the special case of fiber squares, we can give a set of
linearly independent elements in the N\'eron-Severi groups;  moreover, if
$
\E\rarr C
$
is an elliptic modular surface these elements form a basis under Tate's
conjecture.  This is the fiber square analog of the Shioda-Tate Theorem.
}
\fi

\vsp

Our first result gives an analog of the analytic form of Nagao's conjecture
for higher genus fibrations  (see Section \ref{sec:raynaud} for the definition
of Chow traces, and see Section \ref{sec:nagao} for the definition of
$
a_\p(S_t)
$).
Note that the hypothesis on multiplicities in Theorem \ref{thm:jac} is
satisfied in the case of a semistable fibration, or if the fibration has a
section.

\vsp

\begin{thm}
	\label{thm:jac}
Let $S$ be a smooth projective surface, and let
$
\pi: S\rarr C
$
be a fibration over a smooth projective curve $C$, all defined over a number
field $k$.  Suppose that for every point
$
x\in C(\ov{k})
$,
the  GCD of the multiplicities of the irreducible components of the fiber
$
\pi^{-1}(x)
$
is $1$, and that the Chow trace of this fibration is zero.  Assume Tate's
conjecture for $S/k$.  Then
\begin{eqnarray*}
\raisebox{-7pt}{$\stackrel{\displaystyle\text{res}}{\scriptstyle s=1}$}\,
\sum_\p
\frac{-\log N_\p}{N_\p^{s+1}}\!\!\sum_{t\in C(\F_\p)}a_\p(S_t)
&=&
	\begin{array}{ll}
	\text{the rank of the group of sections of the}
	\\
	\text{Jacobian fibration associated to $S\rarr C$.}
	\end{array}
\end{eqnarray*}
\end{thm}

\begin{remark}
For any integer $g\ge 1$ there exists a semistable, genus $g$ fibration
$
X\rarr \mathbf P^1
$
over some number field, such that $X$ is a rational surface and that the
Mordell-Weil rank of the Jacobian fibration is $4g+4$;
	cf.~\cite{saito} and \cite{shioda:manin}.
Since Tate's conjecture is true for rational surfaces, this shows that
Theorem 
	\ref{thm:jac}
is not vacuous.
%

\end{remark}

\vsp

Next, consider the problem of studying Nagao sums for
$
a_\p(\cale_t)^n
$
where $\E\rarr C$ is an elliptic fibration.  Since
$
|a_\p(\cale_t)^n| \le (4N_\p)^{n/2}
$,
\begin{eqnarray}
\underset{s=1}{\text{ res }}
  \sum_\p
    \frac{-\log N_\p}{N_\p^{s+\lambda}}   \sum_{t\in C(\F_\p)} a_\p(\E_t)^n
= 0
&&
	\text{for every real number $\lambda>n/2+1$.}
\label{trivial}
\end{eqnarray}
It is then natural to try to compute the residue at
$
\lambda = n/2+1
$,
and to ask if there is any geometric interpretations of these residues
at every integer
$
\ge n/2+1
$.
To state our results, first note that the relation
$
a_\p(\E_t) = p+1-\#\E_t(\F_\p)
$
suggests that the Nagao sum for $a_\p(\E_t)^n$ should be related to
$$
{\mathcal E}^{(n)}
:=
\text{the fiber product of $n$ copies of $\E$ with itself over $C$.}
$$
Indeed, let
$
\efc
$
be a semistable elliptic fibration over a number field $k$.  Deligne
(\cite{deligne}; cf.~also \cite{scholl}) constructs a desingularization
$
\tilde{\mathcal E}^{(n)}
$
of ${\mathcal E}^{(n)}$ via a canonical sequence of blowups.  Set
$
b_j :=$ the number of $k$-rational exceptional divisors of $\ejtilde$.
Since
$
\widetilde{\mathcal E}^{(1)} = \mathcal E
$
is smooth, we have $b_1=0$.  For $\etwo$ the singular locus consists of
finitely many ordinary double points.

\vsp

\begin{thm}
	\label{thm:e2}
Let
$
f: \E\rarr C
$
be a semistable elliptic fibration over $k$ with a section.
Fix an integer $n\ge 2$.  Assume Tate's conjecture for
$
\ejtilde
$
for every
$
1\le j\le n
$.
Then
\begin{eqnarray}
\underset{s=1}{\text{ res }}
  \sum_\p
    \frac{-\log N_\p}{N_\p^{s+n}}   \sum_{t\in C(\F_\p)} a_\p(\E_t)^n
&
=
&
-1
+
\sum_{j=1}^n
	(-1)^j \binomial{n}{j} \cdot
	(
		b_j - \text{rank } NS(\ejtilde/k)
	).
\label{rightside}
\end{eqnarray}
\end{thm}

\begin{remark}
Note that the left side of (\ref{rightside}) makes no reference to any
desingularization of
$
\mathcal E^{(n)}
$,
and that in general 
$
\mathcal E^{(n)}
$
has many desingularizations.  For example, every the ordinary double point on
$
\etwo/\ov{k}
$
admits a small resolution  \cite[Example IV-29]{eisenbud-harris}.
This process replaces each ordinary double point by a
$
\mathbf P^1
$
instead of
$
\mathbf P^1\times\mathbf P^1
$
(as in the case of blowup)
and is an isomorphism outside the double points.
In particular, it does not affect the N\'eron-Severi group.  On the other hand,
when we evaluate the residue
(\ref{rightside})
unconditionally in the case of elliptic modular
surfaces (Theorem
\ref{thm:modular}), we need to work with Deligne's canonical desingularization
in order to relate the threefold
$
\widetilde{\mathcal E_M}^{(2)}
$
to modular forms.

\end{remark}

\vsp

Let us examine the case $n=2$ of this Theorem in greater details.  Since
$
a_\p(\E_t)^2 \le 4N_\p
$,
when $n=2$ the residue in (\ref{rightside}) is between $-4$ and $0$.  Theorem 
\ref{thm:e2} then suggests that the rank of 
$
NS(\etwotilde/k)
$
is equal to
$$
2 \text{ rank } NS(\E/k)
+
(\#\text{$k$-rational singular points of $\etwo$})
-
1 - S_\E,
$$
where
$
S_\E\in \{0, -1, -2, -3, -4\}
$.
On the other hand, there is a natural collection of divisors on
$
\etwotilde
$
obtained by essentially pulling back to $\etwotilde$, via the two projections
$
\etwo\rarr \E
$,
the generators of
$
NS(\E/k)\otimes \Q_l
$
furnished by the Shioda-Tate Theorem; cf.~the list (\ref{list2}) in
	$\S \ref{sec:generator}$
for details.
\if 3\
{
Denote by 
first, viewing $\etwo$ as a subset of
$
\E\times\E
$,
we can write an element of $\etwo$ as a triple $(a, b, x)$ with $a, b\in\E$
such
that $f(a) = f(b) = x$.  By abuse of notation we will think of such a triple
as an element of $\etwotilde$ as well.   Set
$
\mathcal E^{(2)}_\ast := \etwo - \text{$k$-rational singular points}
$.
Consider the following collection of $k$-rational divisors on $\etwotilde$:

\newcounter{property}
\begin{list}
	{{ \rm (\roman{property})}}{\usecounter{property}
				\setlength{\labelwidth}{50pt}}
\item
a `multifiber'
$
\pi^{-1}(\delta)
$,
where $\delta$ is any non-zero $k$-rational divisor of $C$ which is

\noindent
disjoint from the singular locus of $f$;
\item
the blowup divisors
$
\{\beta_s \}_{s\in\Sigma}
$,
where $\Sigma$ denotes the set of $k$-rational singular points of $\etwo$;
\item
the sectional divisors
$
\sigma_{1, i}$  and  $\sigma_{2, j}
$,
which are the Zariski closure in $\etwotilde$ of the $k$-rational
quasi-projective subvarieties
\begin{eqnarray*}
&&\{ (s_i(x), e, x)\in \tilde{\mathcal E}_\ast^{(2)}(\ov{k}):
	e\in\mathcal E(\ov{k}), x\in C(\ov{k})
	\text{ with } f(e)=x
\},
\\
&&
\{ (e, s_j(x), x)\in \tilde{\mathcal E}_\ast^{(2)}(\ov{k}):
	e\in\mathcal E(\ov{k}), x\in C(\ov{k})
	\text{ with } f(e)=x
\},
\end{eqnarray*}
where
$
\{s_i\}_{i\in I}
$
runs through a complete set of independent generators of the group of
\\
sections  of $\efc$ modulo torsions; and
\item
the diagonal divisor $\Delta$, which is the closure in $\etwotilde$ of
$$
\{ (e, e, x)\in \etwostar(\ov{k}): e\in \E(\ov{k}), x\in C(\ov{k})
	\text{ with } f(e)=x\}.
$$
\end{list}
}
\fi
The cardinality of this natural collection of divisors is (cf.~(\ref{number}))
$$
2 \text{ rank }NS(\E/k)  + (\#\text{$k$-rational singular points of $\etwo$}).
$$
In light of Theorem \ref{thm:e2}, if these divisors form a basis of
$
NS(\etwotilde/k)\otimes\Q
$,
then the residue in (\ref{rightside}) would be exactly $-1$ when $n=2$.  We
have the following partial result.

\vsp

\begin{thm}
	\label{thm:divisor}
Let
$
f: \E\rarr C
$
be a semistable elliptic fibration over $k$ with a section.  Then the divisors
in {\rm(\ref{list2})} on $\etwotilde$ give rise to independent elements
in
$
NS(\etwotilde/k)\otimes\Q
$.
\end{thm}

%
%

\begin{cor}
Let
$
f: \E\rarr C
$
be a semistable elliptic fibration over $k$ with a section.  Assume Tate's
conjecture for $\E$ and
$
\etwotilde
$.
Then 
\begin{equation}
\underset{s=1}{\text{ res }}
  \sum_\p
    \frac{-\log N_\p}{N_\p^{s+2}}   \sum_{t\in C(\F_\p)} a_\p(\E_t)^2
<0.
	\label{cor}
\end{equation}
\end{cor}

\begin{remark}
	\label{relation}
Let $\efc$ be a semistable fibration over $k$ with a section.  Suppose that
every $\ov{k}$-irreducible component of the bad fibers is defined over $k$.
Let $\Gamma$ be such a component, and set
$
\gamma := f(\Gamma) \in C
$,
so the fiber $f^{-1}(\gamma)$ is a $m$-gon for some $m\ge 1$.  The fiber above
$\gamma$ in the singular threefold 
$
\etwo \rarr C
$
then consists of $m^2$ divisors.  After blowing up the $m^2$ ordinary
double points, we pick up $m^2$ exceptional divisors along with the original
$m^2$ divisors on
$
\etwotilde
$,
all Cartier.\footnote{Denote
	by
	$
	\mathring{\E}^{(2)}
	$
	the complement of the ordinary double points in	$\etwo$.  Then the
	$m^2$ Weil divisors on $\etwo$ also give rise to $m^2$ divisors
	$\{ W_i\}_i$ on
	$
	\mathring{\E}^{(2)}
	$.
	These are now Cartier divisors since
	$
	\mathring{\E}^{(2)}
	$
	is smooth.  The Zariski closure in $\etwotilde$ of the $W_i$ then
	give rise to $m^2$ Cartier divisors on $\etwotilde$.  
}
%
%
In light of Theorems \ref{thm:divisor} and \ref{thm:e2}, we see that Tate's
conjecture implies that only about half of these divisors are independent
in
$
NS(\etwotilde/k)
$.
What are these  divisorial relations?
\end{remark}

\vsp

The residue in (\ref{cor}) encodes in analytic terms an equidistribution
statement for the traces of Frobenius for the fibers of an elliptic surface:
  To say that the residue is zero is to say
that most of the fibers have small $a_\p$, while to say that the residue is
$-4$ is to say that most of the fibers have $a_\p$ close to the extremal
values
$
\pm 2\sqrt{N_\p}
$.
We can pose the same problem for elliptic curves over number fields; as
usual we get the
strongest results in the context of modular elliptic curves over $\Q$, for
which analytic techniques are available.
Similarly, in the important case of elliptic modular surfaces, by relating the
$a_p(\E_t)^n$-sum to the Selberg trace formula  we prove the following
result.

\vsp

\begin{thm}
	\label{thm:modular}
Fix an odd, positive integer $M$.  Denote by
$
f_M: \E\rarr X_0(M)
$
the elliptic modular surface over $\Q$ associated to the modular curve
$X_0(M)$. 
For every integer $n\ge 2$, denote by $[n/2]$ the largest integer $\le n/2$.
Then the series
$$
\sum_p
\frac{-\log p}
	{p^{s+[n/2]+1}}\!\!\sum_{t\in X_0(M)(\F_p)}
a_p(\mathcal E_t)^{n}
$$
has a meromorphic continuation to the half plane $Re(s)>0$ with at most a
simple pole at $s=1$.  The pole occurs precisely when $n$ is even, in which
case the residue is
$
\displaystyle \frac{-(n!)}{(n/2)!(n/2+1)!}
$.
\end{thm}

\begin{remark}
(a) Since
$
|a_p(\E_t)^n| \le 2^n p^{n/2}
$,
a trivial lower bound of the residue in Theorem \ref{thm:modular} is
$
\ge  -2^n
$.
But the Stirling formula gives
$$
\frac{n!}{(n/2)!(n/2+1)!}
=
\frac{2^n}{(n/2)^{3/2}} \frac{ e^{O(1/n)} }{e^4 \sqrt{2\pi}},
$$ 
which is $o(2^n)$.  Thus on average,
$
|a_p(\E_t)|
$
is `small'.

(b) For odd $n$,  Theorem \ref{thm:modular} is stronger than 
(\ref{trivial}) in that it gives the residue at a point to the \textit{left}
of the trivial boundary of convergence.
Furthermore, for odd $n$, our argument in fact holds over all number fields.
However,
for even $n$ our use of the Selberg trace formula forces us to work over $\Q$.
If we can directly relate, say,
expression (\ref{almost}) to counting points on Kuga fiber varieties without
using the trace formula then we should be able to work with arbitrary base
fields.

(c) The requirement that $M$ be odd is a simplifying assumption;
with additional (tedious) work the Theorem should hold for general $M$.
\end{remark}

%
%
%

\vsp

As a by-product of Theorem \ref{thm:modular}, we get an analytic proof of
Tate's conjecture for the universal elliptic curve $\E_M\rarr X_0(M)/\Q$
and for its fiber square.

\begin{thm}
	\label{thm:tate}
Let $M$ be an odd integer.

{\rm (i)}
The $L$-functions
$
L_2(\E_M/\Q, s)
$
and
$
L_2(\widetilde{\E_M}^{(2)}/\Q, s)
$
associated to the elliptic modular surface $\E_M\rarr X_0(M)$ over $\Q$
has an analytic continuation to $Re(s)> 7/4$ except for a pole at $s=2$ of
order equal to
$
\text{rank} NS(\E_M/\Q)
$
and
$
\text{rank} NS(\widetilde{\E_M}^{(2)}/\Q)
$,
respectively.  In addition, Tate's conjecture holds for
$
\E_M/\Q
$
and for
$
\widetilde{\E_M}^{(2)}/\Q
$.

{\rm (ii)} The collection of divisors in {\rm(\ref{list2})} gives rise to a
basis of
$
NS(\widetilde{\E_M}^{(2)}/\Q)\otimes\Q
$.

\end{thm}

\vsp

\begin{remark}
(a) 
For Kuga fiber varieties over the full congruence subgroup $\Gamma(N)$,
the first part of Tate's conjecture 
	(\ref{conj1})
 has been verified \cite{gordon}; the
argument there is most likely to be applicable for $\E_M$ as well.  The
second part of Tate's conjecture
	(\ref{conj2})
is probably known to the experts for such Kuga varieties too
	(cf.~Remark \ref{remark:geom} in Section \ref{sec:tate}),
but we have not been able to locate a proof in the literature.
%
%
%
%

(b) Theorem \ref{thm:tate}(ii) is the fiber-square analog of the Shioda-Tate
Theorem for
$
\E_M\rarr X_0(M)/\Q
$;
cf.~Theorem \ref{thm:shioda-tate} below.
\end{remark}

\begin{remark}
Our method can be applied to study higher fiber powers of
$
\E_M\rarr X_0(M)
$.
To do that we need an analog of Theorem \ref{thm:divisor} for such higher
fiber powers.  That calls for a
careful accounting of the pull-back divisors obtained via various
projections
$
\ejtilde \rarr \E
$,
together with a detail analysis of the exceptional divisors arised
from Deligne's canonical desingularization.
Such an analysis is also needed if we want to apply our techniques to study
the groups of $d$-cycles on $\ejtilde$  modulo $l$-adic homological
equivalence for $0<d<n$ (for $d=n-1$ we recover the N\'eron-Severi group).
We will address these issues in
a separate paper.
\end{remark}

\begin{remark}
Our method yields no information regarding the torsion subgroups of the
N\'eron-Severi groups of fibered varieties.  For example, are there any
additional torsion divisors on
$
\ejtilde
$
other than the pull-back of torsion divisors on $\E$?  Do non-trivial torsion
divisors on $\E$ pull-back to non-trivial torsion divisors on $\ejtilde$?
\end{remark}

\begin{remark}
	\label{remark:fermat}
Using the computer algebra package PARI, we evaluate numerically the
$a_p^2$-Nagao residues for
several semistable elliptic fibrations  which are not Kuga varieties.
In each case our data strongly suggests that the residue should be $-1$.
In light of the spectral sequence computation in Section \ref{sec:tate},
to settle this question for a general semistable elliptic fibration
$
\efc
$
we need to understanding the $L$-function attached to the $l$-adic
monodromy representation
$
H_\et^1(C/\ov{k}, R^1 f_\ast \Q_l)
$.
In general this is very difficult problem.   Currently, our effort is 
focused on the semistable elliptic surface
$
y^2 = x(x-A^{2n})(x+B^{2n})
$
associated to the Fermat curve
$
A^{2n} + B^{2n} = C^{2n}
$.
This turns out to be the `universal elliptic curve' over the Fermat curve
	\cite{shioda:frey},
the latter being viewed as a modular curve with respect to a non-congruence
subgroup of $SL(2, \Z)$.
We expect that the complex multiplication structure on the Fermat curve will
greatly facilitate the monodromy calculation.

\end{remark}

\begin{remark}
Nagao sums can be formulated for varieties fibered over a base of dimension
$>1$ as well.  To mimic our argument to this setup, among other things we need
to relate the cohomology of the base to that of the fibered variety; cf.~Lemma
\ref{lem:schoen}.  Significant results along this line have been obtained
in the Ph.~D.~dissertation of Wazir \cite{wazir}.
\end{remark}

The discussion above focuses primarily on semistable fibrations.  At the other
end of the spectrum we have the isotrivial, non-split elliptic fibrations:
these are non-split elliptic surfaces
$
\E\rarr C
$
where the $j$-invariant of the generic fiber is an element of $k$.

\vsp

\begin{thm}
	\label{thm:isotrivial}
Let
$
f: \E\rarr C
$
be an isotrivial, non-split elliptic surface over $\Q$.  Then
$$
\raisebox{-7pt}{$\stackrel{\displaystyle\text{res}}{\scriptstyle s=1}$}\,
\sum_p
\frac{-\log p}{p^{s+2}}\!\!\sum_{t\in C(\F_p)}a_p(\mathcal E_t)^2
=
\biggl\{
  \begin{array}{rll}
  0	&&	\text{if $\E\rarr C$ is a cubic, quartic}
  \\
	&&	\text{or sextic twist family,}
  \\
  -1	&&	\text{otherwise.}
  \end{array}
$$
\end{thm}

\begin{remark}
Is there a geometric interpretation of this residue,
\textit{\`a la}
Theorems \ref{thm:divisor} and \ref{thm:tate}(b)?
\end{remark}

\if 3\
{
{\rm (b)} The collection of divisors {\rm(i)} -- {\rm(iv)} are independent
??? generators of
$
NS(\etwotilde/\Q)\otimes\Q
$.

\vsp

\framebox{Tate conjecture for isotrivial fibrations? symmetric squares of
	these?}

}
\fi

\if 3\
{
\vsp

Finally, we connect our results on higher genus fibrations with those on
fiber powers.   Note that Part (a) below can be thought of as an elliptic
surface analog of the Rankin-Selberg convolution.

\begin{thm}
	\label{thm:square}
Let
$
S\rarr C
$
be a genus two Raynaud fibration over $k$ with
$
Tr_{K_C/k}(Jac_{K_C}(\Gamma_\pi)) = 0
$.
Suppose that the associated Jacobian fibration is $k$-isogenous to
a product of two non-split elliptic fibrations
$
\E_i\rarr C
$
over $k$.  Then

{\rm (a)}
$
	\displaystyle
\underset{s=1}{\text{ res }}
  \sum_\p
    \frac{-\log N_\p}{N_\p^{s+2}}   \sum_{t\in C(\F_\p)} a_\p(S_t)^2
=
2\underset{s=1}{\text{ res }}
  \sum_\p
    \frac{-\log N_\p}{N_\p^{s+2}}   \sum_{t\in C(\F_\p)} a_\p(\E_{1, t})a_\p(\E_{2, t}).
$

{\rm (b)}  Suppose further that $\E_1 = \E_2$.  Assume the Tate conjecture
for
$
\E_1, \E_1^{(2)}, S
$
and 
$
S^{(2)}
$; then
\begin{eqnarray*}
\text{rank }NS( \etwotilde/k )
&
=
&
{\textstyle\frac{1}{2}}
\text{rank }NS( \stwotilde/k )
+
(\#\text{singular points of $\etwo$})
\\
&&
+\:
{\textstyle\frac{1}{2}}
-
{\textstyle\frac{1}{2}} (\#\text{singular points of $\stwo$}).
\end{eqnarray*}
\end{thm}

\begin{remark}
Note that the two sides  in Part (b) involve geometric data only.
Can we prove this equality directly?

The Jacobian variety of a genus $2$ curve is birational to the symmetric
power of the curve with itself.   This gives us a
degree $2$ birational map from  $\etwotilde$ to the Jacobian fibration of
$S\rarr C$, which is every suggestive in light of the many factors of $1/2$
in
(\ref{onemore}).  To make this rigorous we need to examine closely the
bad fibers of the N\'eron models of the Jacobian of $S\rarr C$ versus those
of $\etwotilde$.
\end{remark}

\begin{remark}
Given any field $F$ of characteristic $0$, we can find infinitely many
examples of genus two curves over $F$ whose Jacobians are $F$-isogenous
to the product of two $F$-elliptic curves
\cite[p.~154]{cassels}.  Furthermore, we can arrange that
these two elliptic curves to be isomorphic (resp.~not isomorphic) over $F$,
and if $F$ is a function field we can arrange that the elliptic curves be
non-constant,
so the discussion above is not bogus.

	\framebox{explicit construction????}
\end{remark}
}
\fi

\if 3\
{
\begin{notat}
Denote by $\oo_k$ the ring of integers of a number field $k$.  For any maximal
$\p\subset\oo_k$, denote by $\F_\p$ its residue field, and by $N(\p)$ the
absolute norm of $\p$.
Fix an algebraic closure $\ov{k}$ of $k$ and write $G_k = \gal(\ov{k}/k)$.

\end{notat}
}
\fi

\vsp

\subsection{Summary}
Here is a summary of the paper.  In Section \ref{sec:raynaud} we recall
Raynaud's exact sequence and the Shioda-Tate Theorem for higher genus
fibrations.  These geometric data are applied in Sections
	\ref{sec:nagao} and
\ref{sec:power} to derive formulae for the rank of various N\'eron-Severi
groups; Tate's conjecture comes in when we rewrite the Nagao sums in terms of
the logarithmic derivatives of $L$-functions.
In Section \ref{sec:generator} we produce a natural collection of divisors
on 
$
\etwotilde
$;
to prove that they are independent in 
$
NS(\etwotilde/k)\otimes\Q
$
we intersect these divisors on the threefold $\etwotilde$ with various
well-chosen surfaces, reducing the problem of checking independence to the
case of surfaces.  
In Sections \ref{sec:cyclic} and \ref{sec:modular} we switch gears and count
elliptic curves over finite fields with a given level structure and a fixed
$j$-invariant.  This allows us to rewrite, in the case of the universal
elliptic curve
$
\widetilde{\E_M}
$
over
$
X_0(M)/\Q
$,
the $a_p^n$-Nagao sums in terms of the Selberg trace formula, from which
Theorem \ref{thm:tate} follows.  Combined with a spectral sequence
calculation,
this yields an analytic proof of Tate's conjecture for
$
\widetilde{\E_M}^{(2)}
$.
Finally, in Section \ref{sec:isotrivial} we compute the $a_p^2$-Nagao residue
for isotrivial fibrations over $\Q$, by relating these Nagao sums to symmetric
square $L$-functions.

\begin{convent}
We will have many occasions to  reduce a fibered variety $f: V\rarr C$
modulo a prime ideal $\p$.  Unless otherwise stated, the reduction
is taken with respect to some fixed but possibly non-canonical model of the
fibration.  This introduces ambiguity for finitely many $\p$, which does not
affect our computation of the residues of various $L$-functions associated
to $V$.  In a similar vein, we adopt the convention of using only those
$\p$ at which $V, C$ and $f$ all have good reduction with respect to some fixed
model; some fibers of course could be singular.  
\end{convent}

\vsp

\section{Raynaud's exact sequence and the Shioda-Tate formula}
	\label{sec:raynaud}

In this section we recall several results about fibered surfaces; see
\cite{shioda:higher} for more details.

Let $C$ be a smooth projective curve over $k$, and let $S$ be a smooth
projective surface over $k$.  Denote by
$
Jac_k(C)
$
the Jacobian variety of $C/k$, and by
$
\text{PicVar}_{s}(S)
$
the Picard variety of $S/k$; both of these are Abelian varieties over $k$.
Let
$
\pi: S\rarr C
$
be a proper morphism over $k$.  Assume that the generic fibration of
$\pi$ is a smooth projective curve
$
\Gamma_\pi
$
over the function field $k(C)$.  This gives rise to a short exact sequence
$$
\begin{CD}
0 @>>> \text{Pic}^0_{C/k} @>{\pi^\ast}>> \text{Pic}^0_{S/k}
	@>>>
\text{Pic}^0_{S/C},
\\
&&  @|	@|
\\
&& \text{Jac}_k(C)   &&		 \text{PicVar}_k(S)
\end{CD}
$$
where $\text{Pic}_{S/C}$ denotes the relative Picard scheme of $S/C$
	\cite{bosch}.
To identify the image of $\pi^\ast$ we need the notion of Chow trace:
let $A$ be an Abelian variety over the function field
$
K_C := k(C)
$.
The \textit{Chow trace} of $A/K_C$, denoted by
$
Tr_{K_C/k}(A)
$,
is an Abelian variety over $k$ together with a $K_C$-homomorphism
$
\tau: Tr_{K_C/k}(A) \rarr A
$, 
such that for any Abelian variety $B/k$ and a $K_C$-homomorphism
$
\beta: B\rarr A
$,
there exists a unique $k$-homomorphism
$
\beta': Tr_{K_C/k}(A) \rarr B
$
such that
$
\beta = \tau \circ \beta'
$
(cf.~\cite[p.~138]{lang}).
Note that if $f: \E\rarr C$ is a genus one fibration, then $f$ is non-split if
and only if the Chow trace of $\E\rarr C$ is zero.

We say that
$
\pi: S \rarr C
$
is a \textit{Raynaud fibration} if for every fiber of $\pi$ over $\ov{k}$,
the GCD of the
multiplicities of the irreducible components is $1$.  Two notable examples of
Raynaud fibrations are (a) fibrations $S\rarr C$ with a section, and (b)
semistable fibrations.

\vsp

\begin{thm}[\protect{Raynaud \cite[Thm.~2]{shioda:higher}}]
	\label{thm:raynaud}
Suppose
$
\pi: S\rarr C
$
is a Raynaud fibration over $k$.  Then we have a short exact sequence of
Abelian varieties over $k$
\begin{equation}
0 \rarr Jac_{k}(C) \stackrel{\pi^\ast}{\rarr} \text{PicVar}_{k}(S)
	\rarr
Tr_{K_C/k}(Jac_{K_C}(\Gamma_\pi)) \rarr 0.
	\label{seq}
\end{equation}
\end{thm}

\vsp

\begin{cor}
Suppose
$
\pi: S\rarr C
$
is a Raynaud fibration over $k$ with zero Chow trace.  Then for all but
finitely many $\p$,
\begin{equation}
\text{trace} (\text{Frob}_\p, H^i_\et(S\otimes\ov{k}, \Q_l))
=
\text{trace} (\text{Frob}_\p, H^i_\et(C\otimes\ov{k}, \Q_l)).
	\label{apc}
\end{equation}
\end{cor}

Fix a smooth curve $\sigma/k$ on $S$ such that $\pi(\sigma)=C$.  Fix a fiber
$
F/\ov{k}
$
of $\pi$.  Denote by $T_S$ the subspace of
$
NS(S/\ov{k})\otimes \Q
$
generated by $\sigma$ together with all components of all fibers of $\pi$
over $\ov{k}$.  While $\sigma$ need not be a section of $\pi$,  the image in
$
NS(S/\ov{k})\otimes\Q
$ 
of any two choices of $\sigma$ differs by a $\Q$-multiple.  In particular,
$T_S$ is well-defined.

\begin{thm}[\protect{Shioda-Tate \cite{shioda:higher}}]
	\label{thm:shioda-tate}
Suppose
$
\pi: S\rarr C
$
is a Raynaud fibration over $k$ with zero Chow trace.
Then there is a decomposition of $\Q[G_k]$-modules
\begin{equation}
NS(S/\ov{k}) \otimes \Q
	\simeq
\bigl(Jac_{K_C\ov{k}}(\Gamma_\pi)\otimes\Q\bigr) \oplus T_S.
	\label{shioda-tate-isom}
\end{equation}
\end{thm}

\if 3\
{
This readily yields the generalization of (\ref{isom}) we seek.  Raynaud's
condition is weaker than the condition that $\pi$ has a section;
under this stronger condition
Shioda \cite[Thm.~1]{shioda:higher} deduces from (\ref{seq}) the analogue of 
(\ref{shioda-tate}) for higher genus fibrations.  We expect that Shioda's
argument can be modified to work under the weaker hypothesis of Raynaud,
in which case we will obtain (a) and (b) above under Raynaud's condition.
}
\fi

\vsp

\section{Nagao's conjecture for higher genus fibrations}
	\label{sec:nagao}

\begin{proof}[Proof of Theorem \ref{thm:jac}]
For any reduced, irreducible projective curve $X/\F_\p$, not necessarily
smooth, define
$$
a_\p(X) := N_\p + 1 - \#X(\F_\p).
$$
\if 3\
{
We say that
$\pi$
is a \textit{nice fibration} if for all but finitely many $\p\subset\ok$,
if
$
t\in C(\F_\p)
$
then every $\F_\p$-component of the fiber
$$
S_t := \pi^{-1}(t)
$$
is geometrically irreducible and has multiplicity one.
}
\fi
Let $\pi: S\rarr C$ be a fibration of smooth projective surface over a smooth
projective curve, all defined over a number field $k$.  Suppose the fibration
has good
	reduction\footnote{recall our convention}
at $\p$.
For any $t\in C(\F_\p)$, write $S_t$ for the fiber of $\pi$ at $t$.  Define
$$
m(S_t/\F_\p)
:=
\mbox{number of reduced $\F_\p$-rational components of the fiber $S_t$.}
$$
Denote by $\{C_i\}_i$ the $\F_\p$-components of the fiber $S_t$.
Define
$$
a_\p(S_t) 
:=
\sum_i a_\p(C_i)
+
\sum_{x\in S_t(\F_\p)}
	\bigl[\text{(number of $C_i$ passing through $x$)} - 1\bigr]
-
\bigl[ m(S_t/\F_\p)-1 \bigr].
$$
Note that $a_\p(S_t)$ coincides with the usual definition when $S$ is a
genus one fibration
(this
follows from the case-by-case analysis in \cite[Lem.~17]{mike_joe}),
or
when $S_t$ is a smooth fiber (in which case $a_\p(S_t)$ is the
trace of the geometric Frobenius at $\p$ acting on
$
H^1_\et(S_t\otimes\ov{k}, \Q_l)
$).
Moreover, for every
$
t\in C(\F_\p)
$,
we have\footnote{our
	definition of $a_\p$ makes sense even if $S_t/\F_\p$ is
	irreducible but not geometrically irreducible.  For example, if
	$
	S_t/\F_\p
	$
	is the union of two conjugate elliptic curves defined over
	a quadratic extension of $\F_\p$, then
	$
	a_\p(S_t) = N_\p+1
	$,
	so (\ref{st}) gives $\#S_t(\F_\p) = 0$, as it should
	be.}
\begin{equation}
\#S_t(\F_\p) = 1 - a_\p(S_t) + N_\p + (m(S_t/\F_\p)-1) N_\p,
	\label{st}
\end{equation}
which is the analog for our fibration of \cite[Lem.~17]{mike_joe}.
Consequently,
\begin{equation}
\#S(\F_\p)
=
(1+N_\q)\#C(\F_\p)
+
\sum_{t\in C(\F_\p)} a_\p(S_t)
+
\sum_{t\in C(\F_\p)} (m(S_t/\F_\p)-1) N_\p.
	\label{sumovert}
\end{equation}
On the other hand, since $S$ has good
	reduction
at $\p$,
the Weil conjecture gives
\begin{equation}
\#S(\F_\p)
=
\sum_{i=0}^4 (-1)^i \text{trace}(\text{Frob}_\p, H^i_\et(S\otimes\ov{k}, \Q_l)).
	\label{weil}
\end{equation}
Combine (\ref{sumovert}) and (\ref{weil}) together with (\ref{apc}), we get
$$
\sum_{t\in C(\F_\p)} a_\p(S_t)
=
-
\text{trace}(\text{Frob}_\p, H^2_\et(S\otimes\ov{k}, \Q_l))
+
2N_\q
+
 \sum_{t\in C(\F_\p)} (m(S_t/\F_\p)-1)N_\p
$$
Recall that the submodule $T_S$ of the $G_k$-module
$
NS(S/\ov{k})\otimes\Q
$
is generated by a fixed multisection $\sigma/k$ on $S$ along with all
components of all fibers of $\pi$ over $\ov{k}$.  The geometric Frobenius
$
\text{Frob}_\p
$
fixes $\sigma$ and acts non-trivially on any non-$\F_\p$-rational
components of each fiber of $\pi$.  Furthermore, in light of the Weil
conjecture, $\pi$ has at least one smooth fiber over $\F_\p$ for all but
finitely many $\p$, and the image in
$
NS(S/\ov{\F_\p})\otimes\Q
$
of any two such multisections differ by a non-zero $\Q$-multiple. Consequently,
for all but finitely many $\p$,
$$
\sum_{t\in C(\F_\p)} (m(S_t/\F_\p)-1)N_\p
=
2 + \text{trace}(\text{Frob}_\p, T_{S}),
$$
whence
\begin{eqnarray}
\sum_{t\in C(\F_\p)} a_\p(S_t)
&=&
-
\text{trace}(\text{Frob}_\p, H^2_\et(S\otimes\ov{k}, \Q_l))
+
\text{trace}(\text{Frob}_\p, T_{S})N_\p.
	\label{newsum}
\end{eqnarray}
Denote by
$
L(T_S, s)
$
the $L$-function associated to the $G_k$-module $T_S$.  Following
\cite[p.~52]{mike_joe}, we have, for
$
Re(s)>7/4
$,
\begin{eqnarray}
\frac{d}{ds}\log L(T_S, s)
&=&
\sum_\p -\text{trace}(\text{Frob}_\p, T_{S}) \frac{\log N_\p}{N_\p^{s-1}}
	+ O(1),
	\label{log_ts}
\\
\frac{d}{ds}\log L_2(S/k, s)
&=&
\sum_\p
	-\text{trace}(\text{Frob}_\p, H^2_\et(S\otimes\ov{k}, \Q_l))
	\frac{\log N_\p}{N_\p^{s-1}}
	+ O(1).
	\label{log_l2}
\end{eqnarray}
Combine (\ref{newsum}), (\ref{log_ts}) and (\ref{log_l2}), we get
\begin{equation}
\raisebox{-7pt}{$\stackrel{\displaystyle\text{res}}{\scriptstyle s=1}$}\,
\sum_p
\frac{-\log N_\p}{N_\p^{s+1}}\!\!\sum_{t\in C(\F_\p)}a_\p(\mathcal S_t)
=
-
\raisebox{-7pt}{$\stackrel{\displaystyle\text{ord}}{\scriptstyle s=2}$}\,
	L_2(S/k, s)
+
\raisebox{-7pt}{$\stackrel{\displaystyle\text{ord}}{\scriptstyle s=1}$}\,
	L(T_S, s).
	\label{final}
\end{equation}
$T_S$ is a finite dimensional $\Q[G_k]$-module, so $L(T_S, s)$ is an Artin
$L$-function.  In particular, the second term on the right side of
(\ref{final}) is
$
-\text{rank}T_S(k)
$.
Under Tate's conjecture, the first term on the right side of (\ref{final})
is
$
\text{rank}NS(S/k)
$.
Invoke the Shioda-Tate isomorphism (\ref{shioda-tate-isom}) and we get
Theorem \ref{thm:jac}.

\end{proof}

\if 3\
{
\section{Fiber-square of elliptic fibrations}

\begin{proof}[Proof of Theorem \ref{thm:e2}]
Let $\efc$ be a semistable elliptic fibration over a number field $k$.  Then
the singular locus of
$
\etwo := \E\times_C \E
$
consists of isolated ordinary double points only.  Denote by
$
\etwotilde
$
the desingularization of $\etwo$ by blowing up these double points.
Set
$
b_\p(\E_t) = \text{trace}(\text{Frob}_\p, H^2_{et}(\E_t\otimes\ov{k}, \Q_l))
$
and
$
t_{\p,i}(\etwotilde)
=
\text{trace}(\text{Frob}_\p, H^i_{et}(\etwotilde_t\otimes\ov{k}, \Q_l))
$.
Since
$
a_\p(\E_t)^2 = (N_\p+1 - \#\E_t(\F_\p))^2
$
and since the difference between
$
\sum_{t\in C(\F_\p)} \E_t(\F_\p)^2
$
and
$
\#\E^{(2)}(\F_\p)
$
comes from the singular points only, a little arithmetic gives
\begin{eqnarray}
\sum_{t\in C(\F_\p)} a_\p(\E_t)^2
&=&
N_\p^2+N_\p
-
(N_\p^2 + 2N_\p)
\times
(\#\text{singular points of $\etwo$})
	\label{main}
\\
&&
+
	t_{\p,2}( \etwotilde )
-
	2b_\p(\E)
+
	2N_\p a_\p(C)
-
	t_{\p,3}( \etwotilde )
\label{second}
\\
&&
+
	N_\p\bigl(
	  t_{\p, 2}( \etwotilde )  - 2b_p(\E)
	\bigr),
\label{third}
\end{eqnarray}
Note that every term in (\ref{second}) is
$
\ll_{\E, C} N_\p^{3/2}
$.
Consequently,
\begin{eqnarray}
\lefteqn{
\underset{s=1}{\text{ res }}
  \sum_p
    \frac{-\log N_\p}{N_\p^{s+2}}   \sum_{t\in C(\F_\p)} a_\p(\E_t)^2
\:=\:
-1
-
\bigl(
  \#\text{singular points of $\etwo$}
\bigr)
}
	\nonumber
\\
&&
-
\; 2 \underset{s=1}{\text{ res }}
  \sum_p
    \frac{-\log N_\p}{N_\p^{s+1}} b_\p(\E)
+
\underset{s=1}{\text{ res }}
  \sum_\p
    \frac{-\log N_\p}{N_\p^{s+1}} t_{\p, 2}( \etwotilde ).
\label{res_1}
\end{eqnarray}

By Nagao's conjecture for $\efc$, the first residue in (\ref{res_1}) is
the rank of
$
\E(K_C)
$,
while the Tate conjecture for
$
\etwotilde
$
predicts that the second residue in (\ref{res_1}) is the rank of
$
NS( \etwotilde/k )
$.
This completes the proof of Theorem \ref{thm:modular}.

\end{proof}

}
\fi

\vsp

\section{Fiber powers of semistable elliptic fibrations}
	\label{sec:power}

\begin{proof}[Proof of Theorem \ref{thm:e2}]
Let $\efc$ be a semistable elliptic fibration over a number field $k$.  Fix
an integer $n\ge 2$.  Then
\begin{eqnarray}
\sum_{t\in C(\F_\p)} a_\p(\E_t)^n
&=&
\sum_{t\in C(\F_\p)} (N_\p + 1 - \#\E_t(\F_\p) )^n
	\nonumber
\\
&=&
\sum_{j=0}^n
	(-1)^j \binomial{n}{j} (N_\p+1)^{n-j}
\sum_{t\in C(\F_\p)}
	\#\E_t(\F_\p)^j
	\nonumber
\\
&=&
\#C(\F_\p) (N_\p+1)^n
+
\sum_{j=1}^n
	(-1)^j \binomial{n}{j} (N_\p+1)^{n-j} \cdot \#\ej(\F_\p),
	\label{fix}
\end{eqnarray}
where $\ej$ denotes the fiber product of $j$ copies of $\E$ with itself over
$C$.  Denote by
$
\ejtilde
$
Deligne's desingularization of $\ej$ \cite{deligne} via a canonical
sequence of blowups.  See \cite[$\S 2$]{scholl} for a description of the
desingularization.  Write
$$
b_j := \text{the number of $k$-rational exceptional divisors of $\ejtilde$}.
$$
Since $\E^{(1)}=\E$ is smooth, we have $b_1=0$.  For $\etwotilde$ the singular
locus consists of isolated ordinary double points.
Write
$
O_\E( \: \cdot \:)
$
for the usual big-$O$ notation with the constant depending
on $\E$ only (and not on $\p$); as usual the constant could change from line
to line.  
Then
$$
\#\ej(\F_\p) = \#\ejtilde(\F_\p) + b_j(N_\p^j + O_\E(N_\p^{j-1})).
$$
Substitute this into (\ref{fix}) and expand
$
(N_\p+1)^{n-j}
$
using the binomial theorem, we get
\begin{eqnarray}
&&
(N_\p+1)^n (N_\p + 1 - a_\p(C))
	\nonumber
\\
&&
+
\sum_{j=1}^n
	(-1)^j \binomial{n}{j} \sum_{l=0}^{n-j} \binomial{n-j}{l} N_\p^l
	\times
	\Bigl[
		\#\ejtilde(\F_\p)
		-
		(N_\p^j + O_\E(N_\p^{j-1})) \cdot b_j
	\Bigr].
	\nonumber
\end{eqnarray}
Set
$
tr_\p^i( \ejtilde )
:=
\text{trace}( \text{Frob}_\p, H^i(\ejtilde, \Q_l) )
$.
Then
$
\#\ejtilde(\F_\p) = \sum_{j=0}^{2j+2} (-1)^j tr_\p^i( \ejtilde )
$
by the Weil conjecture, so (\ref{fix}) becomes
\begin{eqnarray}
&&
N_\p^{n+1} + (n+1)N_\p^n - N_\p^n a_\p(C) + O_\E(N_\p^{n-1})
	\nonumber
\\
&&
+
\sum_{j=1}^n
	(-1)^j \binomial{n}{j} \sum_{l=0}^{n-j} \binomial{n-j}{l} N_\p^l
	\times
	\Bigl[
		\sum_{i=0}^{2j+2} (-1)^i tr_p^i (\ejtilde)
		-
		b_j N_\p^j + O_\E(N_\p^{j-1})
	\Bigr]
	\nonumber
\\
&=&
N_\p^{n+1} + (n+1)N_\p^n - N_\p^n a_\p(C) + O_\E(N_\p^{n-1})
	\label{11}
\\
&&
+
\sum_{j=1}^n
	(-1)^j \binomial{n}{j} \sum_{l=0}^{n-j} \binomial{n-j}{l} N_\p^l
	\times
	\sum_{i=0}^{2j+2} (-1)^i tr_p^i (\ejtilde)
	\label{22}
\\
&&
-
N_\p^n \displaystyle \sum_{j=1}^n (-1)^j \binomial{n}{j} b_j + O_\E(N_\p^{n-1}).
	\label{33}
\end{eqnarray}

\begin{lem}
	\label{lem:schoen}
Suppose $f$ has a section.  Then for every $j>0$ we have
$
tr_\p^1(\widetilde{\mathcal E}^{(j)})
=
a_\p(C)
$.
\end{lem}

\begin{proof}
The semistable fibration $\efc$ induces a fibration
$
\tilde{p}(j): \ejtilde \rarr C
$.
Schoen \cite[Lem.~1]{schoen} shows that $\tilde{p}(j)$ induces an
isomorphism of fundamental groups
$
\pi_1(\ejtilde\otimes\C) \rarr \pi_1(C\otimes\C)
$
when $j=2$; his
	argument\footnote{which calls for $f$ to have a section.}
applies
	\textit{mutatis mutandis}
for $j>2$ as well.  By the compatibility of $l$-adic cohomology with
Betti cohomology, 
$
\tilde{p}(j)_\ast: H^1_\et(\ejtilde\otimes\ov{k}, \Q_l) \rarr H^1_\et(C\otimes\ov{k}, \Q_l)
$
is then a $G_k$-equivariant isomorphism for $j>1$.  When $j=1$, Raynaud's
exact sequence (\ref{seq}) implies that
$
\tilde{p}(j)_\ast
$
is also an isomorphism.  The Lemma then follows.

\end{proof}

\if 3\
{
To prove the claim we adopt the argument in \cite[Lem.~1]{schoen}, which
handles the case $j=2$.  We work over the complex numbers.  Denote by
$
C^\circ
$
the complement in $C$ of the singular locus of the semistable fibration
$
\efc
$,
and by
$
f^\circ: \E^\circ \rarr C^\circ
$
the induced smooth fibration.  Denote by
$
{\ej}^\circ
$
the fiber product of $j$ copies of $\E^\circ$ over $C^\circ$, and by
$
p^\circ: {\ej}^\circ \rarr C^\circ
$
the induced smooth fibration.   Fix a point $x\in C^\circ$.  We have a
commutative diagram
$$
\begin{CD}
\pi_1( p^{-1}(x) )  @>{\phi}>>  \pi_1({\ej}^\circ)  @>{p^\circ_\ast}>> \pi_1(C^\circ)
\\
&&	@V{i}VV		@VV{j}V
\\
&&
\pi_1(\ejtilde)  @>{\tilde{p}_\ast}>> \pi_1(C)
\end{CD}
$$
where the top row is exact, both $i, j$ are onto, and
$
\tilde{p}_\ast
$
splits since $\efc$ has a section.  Thus it remains to check that
$
\tilde{p}_\ast
$
is injective.

Suppose $\gamma$ is a loop in $\ejtilde$ whose class $[\gamma]$ in 
$
\pi_1(\ejtilde)
$
goes to $0$ under $\tilde{p}_\ast$.  Then ??????
Thus it remains to show that
$
\pi_1( p^{-1}(x) )
$
has trivial image in 
$
\pi_1(\ejtilde)
$.

Let
$
Y_l \subset \ej - \ej_{\text{sing}}\subset \ej
$
be a section for the canonical projection
$
pr_l: \ej \rarr \E
$.
Then
$$
\pi_1(p^{-1}(x)\cap Y_1) \times \cdots \times \pi_1(p^{-1}(x)\cap Y_j)
\simeq
\pi_1(p^{-1}(x)),
$$
so we are reduced to show that
$
\pi_1(f^{-1}(x))
$
has trivial image in 
$
\pi_1(\E)
$.

}
\fi

Invoke the Lemma and we see that 
$
(\ref{11}) + (\ref{22})
$
is equal to

\begin{eqnarray*}
&&
N_\p^{n+1} + (n+1)N_\p^n - N_\p^n a_\p(C) + O_\E(N_\p^{n-1})
\\
&&
+
\sum_{j=1}^n
	(-1)^j \binomial{n}{j}
	\Big[
		N_\p^{n-j} + (n-j)N_\p^{n-j-1} + O_\E(N_\p^{n-j-2}) 
	\Bigr]
	\times
\\
&&
\hspace{88pt}
	\Bigl[
		N_\p^{j+1}
		-
		N_\p^j a_\p(C)
		+
		N_\p^{j-1} tr_\p^2(\ejtilde)
		+
		O_\E(N_\p^{j-1/2})
	\Bigr]
\\
&=&
N_\p^{n+1} + (n+1)N_\p^n - N_\p^n a_\p(C) + O_\E(N_\p^{n-1})
\\
&&
+
\sum_{j=1}^n
	(-1)^j \binomial{n}{j}
	\Big[
		N_\p^{n+1} - N_\p^n a_\p(C) + (n-j) N_\p^n
		+
		N_\p^{n-1} tr_\p^2(\ejtilde)
		+
		O_\E(N_\p^{n-1/2})
	\Bigr]
\\
&=&
\bigl(
	N_\p^{n+1} - N_\p^n a_\p(C)
\bigr)
\sum_{j=0}^n (-1)^j \binomial{n}{j}
+
N_\p^n
\Bigr[
	n+1+\sum_{j=1}^n (-1)^j \binomial{n}{j}(n-j)
\Bigr]
\\
&&
+
N_\p^{n-1}
\sum_{j=1}^n
	(-1)^j \binomial{n}{j} tr_\p^2(\ejtilde)
+
O_\E(N_\p^{n-1/2})
\\
&=&
N_\p^n 
+
N_\p^{n-1}
\sum_{j=1}^n
	(-1)^j \binomial{n}{j} tr_\p^2(\ejtilde)
+
O_\E(N_\p^{n-1/2})
\end{eqnarray*}
since
$
\sum_{j=0}^n (-1)^j
  \bigl(
    \begin{smallmatrix}n\\j\end{smallmatrix}
  \bigr)
  (n-j)
=0
$
for $n>1$.  Thus
\begin{eqnarray}
\raisebox{-7pt}{$\stackrel{\displaystyle\text{res}}{\scriptstyle s=1}$}\,
	\sum_\p
	\frac{-\log N_\p}{N_\p^{s+n}}
\sum_{t\in C(\F_\p)} a_\p(\E_t)^n
&=&
\raisebox{-7pt}{$\stackrel{\displaystyle\text{res}}{\scriptstyle s=1}$}\,
	\sum_\p
	\frac{-\log N_\p}{N_\p^s}
	\Bigl(
	  1 - \sum_{j=1}^n (-1)^j \binomial{n}{j} b_j
	\Bigr)
		\label{firstterm}
\\
&&
+
\sum_{j=1}^n
	(-1)^j \binomial{n}{j} 
\raisebox{-7pt}{$\stackrel{\displaystyle \text{res}}{\scriptstyle s=1}$}\,
	\sum_\p
	\frac{-\log N_\p}{N_\p^{s+1}}tr_\p^2(\ejtilde)
		\label{secondterm}
\\
&=&
-1
+
\sum_{j=1}^n
	(-1)^j \binomial{n}{j}
	\Bigl[
\raisebox{-7pt}{$\stackrel{\displaystyle \text{ord}}{\scriptstyle s=2}$}\,
	L_2(\ejtilde, s) + b_j
	\Bigr].
		\label{thirdterm}
\end{eqnarray}
Apply Tate's conjecture and we get Theorem \ref{thm:e2}.

\end{proof}

\section{N\'eron-Severi groups of fiber squares}
	\label{sec:generator}

Let $\efc$ be a semistable fibration over $k$.  Fix a non-zero, $k$-rational
divisor
$
d_0
$
on $C/k$ over which $f$ is smooth.  Set
$
F_\E := f^{-1}(d_0)
$.
A different choice of $d_0$ gives rise to another vertical fiber whose image
in
$
NS(\E/k)\otimes \Q
$
is a non-zero $\Q$-multiple of that of $F_\E$.

Denote by
$
\{ t \}_{t\in \tau}
$
the set of $k$-rational, irreducible divisors of $C/k$ over which $\efc$ is
singular.  Since $f$ is semistable, for every $t$ the fiber
$
f^{-1}(t)
$
is a sum
$
\sum_{i=1}^{m(t)} A_i(t)
$
of pairwise distinct, $k$-irreducible divisors on the surface $\E$.  Suppose
further than $\E$ has a zero section $0_\E$.  Then for each
$
t\in\tau
$
exactly one $A_i(t)$ intersects non-trivially with $0_\E$.  Relabel the
components if necessary, we can assume that
$
A_1(t)\cap 0_\E \not=\emptyset
$
for every $t$.  The Shioda-Tate Theorem (Theorem \ref{thm:shioda-tate}) says
that
$
NS(\E/k)\otimes\Q
$
has a basis consisting of
	\newcounter{property}
\begin{list}
	{{ \rm (\roman{property})}}{\usecounter{property}
				\setlength{\labelwidth}{50pt}}
\item
a set of generator $\{ s_i\}_{i\in I}$ of the group of sections of
	$\efc$ modulo torsions,
\item
the zero section $0_\E$,
\item
a vertical fiber $F_\E$, and
\item
the components	$\{ A_i(t) \}_{i>1}$	for each $t\in \tau$.
	\label{component}
\end{list}
Since
$
f^{-1}(t) = \sum_{i=1}^{m(t)} A_i(t)
$
is a non-zero multiple of $F_\E$ in $NS(\E/k)\otimes\Q$, for
	({\roman{property}})
above we could leave out any one of the $A_i(t)$, not just $A_1(t)$; we choose
$A_1(t)$ to facilitate future calculations.

Denote by
$
\pi_i: \etwo\rarr \E
$
the two projections
$$
\begin{CD}
\E^{(2)}  @>{\pi_1}>> \E
\\
@V{\pi_2}VV	@VV{f}V
\\
\E @>{f}>> C.
\end{CD}
$$
Then we have a natural collection of $k$-rational divisors on
$
\etwo
$:
\begin{equation}
	\renewcommand{\arraystretch}{1.2}
\left\{
\begin{array}{lll}
S_{ij} &:=& \pi_i^{-1}(s_j),	\text{ where $i=1, 2$ and $j\in I$,}
\\
A_{ij}(t) &:=& \pi_i^{-1}A_j(t),
		\text{ where $i=1, 2$, and $1<j\le m(t)$ with $t\in \tau$,}
\\
F^{(2)} &:=& F_\E\times_C F_\E,
\\
E_i &:=& \pi_i^{-1}(0_\E),	\text{ $i=1, 2$,}
	\hspace{10pt}
	\text{ and}
\\
\Delta &:=& \{ (e, e)\in \etwo\subset \E\times\E: e\in \E \},
\end{array}
\right.
	\label{list1}
	\renewcommand{\arraystretch}{1}
\end{equation}
where in the definition of $\Delta$ we view $\etwo$ as a subset of
$
\E\times\E
$.
The singular locus of $\etwo/\ov{k}$ consists of finitely many ordinary
double points.  These points break up into a finite set of $k$-conjugate
orbits
$
\{ l \}_{l\in L}
$.
Denote by
$
\beta: \etwotilde\dasharrow \etwo
$
the blowup of $\etwo$ at \textit{all} of these double points, and for each
$l\in L$, denote by $B_l$ the exceptional divisor over $l$.    Theorem
	\ref{thm:divisor}
would follow if we can show that the following collection of divisors
\begin{equation}
	\renewcommand{\arraystretch}{1.2}
\biggl\{
\begin{array}{ll}
\tilde{S}_{ij} := \beta^{-1}(S_{ij}),
\:\:
\tilde{A}_{ij}(t) := \beta^{-1}(A_{ij}(t)),
\:\:
\tilde{F}^{(2)} := \beta^{-1}(F^{(2)}),
\\
\tilde{E}_i := \beta^{-1}(E_i),
\:\:
\tilde{\Delta} := \beta^{-1}(\Delta),
\:\:
\text{and the blowup divisors }  B_l 
\end{array}
\biggr.
	\renewcommand{\arraystretch}{1}
	\label{list2}
\end{equation}
are linearly
	independent\footnote{Instead
	of 
	the $A_{ij}(t)$, it might seem more natural to work with
	$
	A_{j_1}(t)\times_C A_{j_2}(t)
	$.
	However, these elements turn out to be dependent in
	$NS(\etwotilde/k)\otimes\Q$; cf.~Remark \ref{relation} in the
	Introduction.}
in $NS(\etwo/k)\otimes\Q$.  Note that the number of
these divisors is
\begin{eqnarray}
&&
2 (\#I) + \sum_{t\in\tau} 2( m(t) - 1) + 4 + (\#L)
	\nonumber
\\
&=&
2 \text{ rank } NS(\E/k) + (\#\text{$k$-rational points of $\etwo$})
	\label{number}
\end{eqnarray}
by the Shioda-Tate Theorem.  Since blowing up the $k$-Galois orbit of an
ordinary
double point increases by one the rank of the Picard group over $k$, and hence
the N\'eron-Severi group over $k$, it suffices to
show that the divisors in the list (\ref{list1}) are linearly independent in
	$NS(\etwo/k)\otimes\Q$.  Note that these divisors are all Cartier:
$
A_{ij}, S_{ij}, F^{(2)}
$
and
$
E_i
$
are pull-back of Cartier divisors on $\E$, and $\Delta$ is the intersection of
a Cartier divisor on $\E\times\E$ with the subvariety $\etwo$.

Suppose there is a relation of the form
$$
R := \sum_t \sum_{i=1, 2}\sum_{j>1} \alpha_{ij}(t) A_{ij}(t)
+
\sum_{i=1, 2}\sum_{j>1} \sigma_{ij} S_{ij}
+
\gamma F^{(2)}
+
\delta \Delta
+
\sum_{i=1, 2} \epsilon_i E_i
= 0
	\hspace{10pt}
	\text{ in $NS(\etwo)$}
$$
with
$
\alpha_{ij}(t), \sigma_{ij}, \gamma, \delta, \epsilon_i\in \Q
$.
We can find a non-zero, $k$-rational divisor $d_1$ on $C$ which is disjoint
from every
$
t\in\tau
$
and from $d_0$ (recall that $F_\E := f^{-1}(d_0)$), and that
$
G := f^{-1}(d_1)
$
intersects transversely with every $S_{ij}, \Delta$ and $E_i$.  Then
$$
G\cdot A_{ij} = G\cdot F = 0.
$$
The relation $R\cdot G = 0$ then says that
\begin{equation}
\sum_{i=1, 2}\sum_{j>1} \sigma_{ij} (S_{ij}\cdot G)
+
\delta (\Delta\cdot G)
+
\sum_{i=1, 2} \epsilon_i (E_i\cdot G)
=
0
	\label{rg}
\end{equation}
in $NS(G/k)\otimes\Q$.  Since $f$ is semistable, infinitely many of the fibers
of $f$ (over $\ov{k}$) are elliptic curves without complex multiplication.
Thus we can assume that
$
G\simeq \E_{d_1}\times \E_{d_1}
$,
where
$
\E_{d_1} := f^{-1}(d_1)
$
is the Galois orbit over $k$ of an elliptic curve over $\ov{k}$ without
complex multiplication.  Consequently,
$
NS(\E_{d_1}\times \E_{d_1}/k)\otimes\Q
$
is generated by three independent $k$-rational divisors, namely
$
\E_{d_1}\times 0, 0\times \E_{d_1}
$
and the diagonal divisor on
$
\E_{d_1}\times \E_{d_1}
$;
and we have the identification
$$
\begin{array}{llll}
S_{1j}\cdot G \sim 0\times \E_{d_1},
&
E_1\cdot G    \sim \E_{d_1}\times 0,
&
\Delta\cdot G \sim \text{ diagonal divisor on
$
\E_{d_1}\times \E_{d_1}
$,}
\\
S_{2j}\cdot G \sim \E_{d_1}\times 0,
&
E_2\cdot G    \sim 0\times \E_{d_1}.
\\
\end{array}
$$
Then (\ref{rg}) implies that (among other things)
\begin{equation}
\delta = 0.
	\label{relation_one}
\end{equation}
Next, consider the intersection of $R$ with $E_1\simeq \E/k$.  Note that
$$
S_{1j}\cdot E_1 = 0	\hspace{10pt}\text{ for all $j$,}
$$
and we have the identification
$$
\begin{array}{lllll}
A_{2j}(t) \cdot E_1	\sim  A_j(t),
&
S_{2j}\cdot E_1		\sim s_j,
&
F^{(2)}\cdot E_1	\sim F_\E,
&
E_2 \cdot E_1		\sim 0_\E
\end{array}
$$
in $NS(E_1/k)\otimes\Q \simeq NS(\E/k)\otimes\Q$.  Recall that
$
A_j(t)
$
intersects trivially with $0_\E$ for all $j>1$,  so
$$
A_{1j}(t)\cdot E_1 = 0	\hspace{10pt}\text{ for all $j>1$.}
$$
We claim that $E_1\cdot E_1=0$; in light of the identifications above the
relation $R\cdot E_1=0$ then becomes
$$
\sum_t \sum_{j>1} \alpha_{2j}(t) A_j(t)
+
\sum_{j>1} \sigma_{2j} s_j
+
\gamma F_\E
+
\epsilon_2 0_\E
= 0
$$
in $NS(E_1/k)\otimes\Q \simeq NS(\E/k)\otimes\Q$.  The Shioda-Tate Theorem
then implies that
$$
\sigma_{2j} = \alpha_{2j} = \gamma = \epsilon_2 = 0.
$$
for every $j$.  Repeat the same argument for $R\cdot E_2$ and we get that
$
\sigma_{1j} = \alpha_{1j} = \epsilon_1 = 0
$
for every $j$.  All in all, these calculations show that the divisors in
(\ref{list2}) are independent modulo numerical equivalence, and hence they
are independent modulo algebraic equivalence tensored with $\Q$.

It remains to show that $E_1\cdot E_1 = 0$ in $NS(\etwo/k)\otimes\Q$.  For that
we need an auxilary result.

\begin{lem}
	\label{lem:aux}
Let $\E\rarr C$ be an elliptic fibration over $k$ with a zero section
$
0_\E
$.
Then there exists a divisor $M$ on $\E$ which is disjoint from $0_\E$,
such that $M-n0_\E$ is the divisor
of a function on $\E$ for some integer $n>1$.
\end{lem}

\begin{proof}
We can assume that $\E\rarr C$ is relatively minimal.  Write $\E'$ for the
generic fiber.  Fix an integer $m>1$.  Then
\begin{eqnarray*}
M'
&:=&
\text{kernal of the multiplication-by-$m$ map on $\E'$,}
\end{eqnarray*}
and hence
$
M'' := M' - \{0_{\E'}\}
$,
are both divisors on the elliptic curve $\E'$ over the function field $K_C$
of $C/k$.  Add up the points in $M''$ and we see that
$
M'' - (m^2-1)0_{\E'}
$
is the divisor of a function on $\E'/k_C$.
Since $\E\rarr C$ is relatively minimal, by the N\'eron property $M''$ extends
to a divisor $M$ on $\E$ which is disjoint from $0_\E$, and the Lemma follows
by taking $n = m^2-1$.

\end{proof}

Returning to the proof of $E_1\cdot E_1=0$, let $M$ and $n$ be as in the Lemma.
Then $E_1$ and 
$
\E\times_C M
$
defines the same elements in
$
NS(\etwo/k)\otimes\Q
$,
whence
$
E_1\cdot E_1
$
is a non-zero $\Q$-multiple of
$
(\E\times_C 0_\E)\cdot (\E\times_C M)
$,
which is zero.  This completes the proof that the divisors in (\ref{list1})
are linearly independent in
$
NS(\etwo/k)\otimes\Q
$,
and Theorem \ref{thm:divisor} follows.

	\qed

\vsp

\if 3\
{

\framebox{blowup an ordinary double point increases the rank of the Picard
	group by one??}

\vsp

	\noindent
	\framebox{Step 2.}

why do the divisors in (\ref{list1}) remain independent after blowup?
	\framebox{a la Cox?}

\vsp

	\noindent
	\framebox{Step 3.}

\begin{itemize}
\item
why are the blowup divisors independent from these?
\item
triple self-intersection?
\end{itemize}
}
\fi

%
%
%
%
%
%

\if 3\
{
\PSbox{cone01.ps hoffset= 50 voffset=100}{1in}{1.5in}

\PSbox{cone02.ps hoffset=260 voffset=208}{1in}{1.5in}

\newpage

\PSbox{cone03.ps hoffset= 50 voffset=100}{1in}{1.5in}

\PSbox{cone0.ps  hoffset=260 voffset=208}{1in}{1.5in}
}
\fi

\vsp

\section{Cyclic subgroups of elliptic curves}
	\label{sec:cyclic}

The classical modular curve $X_0(M)$ has an open subset
$
Y_0(M)\subset X_0(M)
$,
both defined over $\Q$, such that for any field $K$ of characteristic prime
to $N$, the $K$-rational points of $Y_0(M)$ parameterizes equivalent classes
of pairs
$
(E, C)
$,
where $E$ is a elliptic curve over $K$, and $C$ is a $K$-rational cyclic
subgroup of $E$ of order $M$; two such pairs
$
(E, C)
$
and
$
(E', C')
$
are declared to be equivalent if there exists a $K$-rational isomorphism
taking $E$ to $E'$ and $C$ to $C'$.    See \cite{ogg} for more details.

\if 3\
{
For any elliptic curve $E$ over a field, denote by $j(E)$ its $j$-invariant.
The following result is an immediate consequence of the description of
$X_0(M)$ above
	(cf.~\cite[p.~40]{ito}).

\begin{lem}
	\label{lem:twist}
Suppose $p>3$.  Then the number of pairs $(E, C)$ over $\F_p$ corresponding a
given $\F_p$-rational point on
$
X_0(M)
$
and with $j(E)\not=0, 1728$ is equal to
$
\#{\rm Aut}_{\ov{\F}_p}(E)
$.

	\qed
\end{lem}
}
\fi

There is a canonical map
$
X_0(M) \rarr X_0(1)
$
defined over $\Q$ which sends a pair $(E, C)\in Y_0(M)$ to
$
j(E) :=
$
the $j$-invariant of $E$.  In preparation for the next Section we need to
understand the number of $\fp$-rational points on each fiber of this map.
Ogg's argument in \cite[Thm.~2]{ogg} applies \textit{mutatis mutandis} and
yields the following result.

\begin{lem}
	\label{lem:ogg}
Fix $j\in\F_p$ with $j\not=0, 1728$. Choose any elliptic curve $E/\F_p$ with
$
j(E)=j
$.
Then the number of equivalent classes of pairs $(E', C')$ in
$
Y_0(M)(\F_p)
$
with $j(E')=j$ and $a_p(E)=a_p(E')$, is equal to \underline{one half} of
$$
N_{E, p}(M)
:=
\text{the number of $\F_p$-rational cyclic subgroups of $E$ of order $M$.}
$$
	\qed
\end{lem}

Set $d_{E, p} := a_p(E)^2 - 4p$.  When $M$ is prime to
$
2d_E
$,
Ogg \cite[Prop.~2]{ogg} shows that
\begin{equation}
N_{E, p}(M)
=
\prod_{q|d}\biggl( 1+ \biggl(\frac{d_{E, p}}{q}\biggr) \biggr),
	\label{ogg}
\end{equation}
where $q$ runs through all distinct prime divisors of $M$.   Utilizing
techniques from Waterhouse's thesis \cite{waterhouse}, Ito \cite{ito}
evaluates
$
N_{E, p}(M)
$
in all cases.  For the rest of this section, we state Ito's Theorem and
rewrite his expression for
$
N_{E, p}(M)
$
in a form suitable for our subsequent application.   Clearly it suffices to
work with the case where $M$ is a prime power.  We begin by setting up
some notation.

Let $E/\fp$ be an ordinary elliptic curve, so
$
a := a_p(E)
$
lies between $\pm 2\sqrt{p}$ and  $a\not=0$.  Denotes by
$
\pi_a, \ov{\pi}_a
$
the two distinct roots of $h_{a, p}(x) := x^2-ax+p$.  Fix a prime $l$.  Set
\begin{eqnarray*}
c_a
&:=&
\text{the conductor of the order $\Z[\pi_a]$},
\\
\ea
&:=&
\text{ord}_l(c_a).
\end{eqnarray*}
Then 
$
\text{End}_{\F_p}(E)
$
is an order in the imaginary quadratic field
$
\Q(\pi_a)
$
containing $\Z[\pi_a]$.  The converse also holds.

\begin{lem}[\protect{Waterhouse  \cite[Thm.~4.2]{waterhouse}}]
	\label{lem:waterhouse}
If $a\not=0$, then any order in $\Q(\pi_a)$ containing $\Z[\pi_a]$ is the
endomorphism ring of some elliptic curve over $\F_p$ with
$
p+1-a
$
points over $\F_p$.

	\qed
\end{lem}

\begin{thm}[Ito \cite{ito}]
	\label{thm:ito}
Let $E/\fp$ be ordinary.  With the notation as above, set
\begin{eqnarray*}
c_E
&:=&
\text{the conductor of $\text{End}_{\F_p}(E)$,}
\\
\epsilon_E
&:=&
\text{ord}_l(c_E).
\end{eqnarray*}
Then
$
\epsilon_E\le \ea
$,
and
\begin{list}
	{{\rm(\roman{property})}}{\usecounter{property}
				\setlength{\labelwidth}{30pt}}
\item
if $1\le \epsilon\le \ea-\epsilon_E$ then
$
N_{E, p}(l^\epsilon) = (l+1)l^{\e-1}
$
{\rm;}
\item
if $\ea-\epsilon_E < \e \le \ea+\epsilon_E$, then 
$
N_{E, p}(l^\epsilon) = l^{[(\e+\ea-\epsilon_E)/2]}
$
{\rm;}
\item
If $\e>\ea+\epsilon_E$, then
$$
N_{E, p}(l^\epsilon)
=
\left\{
  \begin{array}{ll}
	2 l^{\ea}	& \text{if $l$ splits in $\Q(\pi_a)$},
	\\
	0	& \text{if $l$ is inert in $\Q(\pi_a)$},
	\\
	l^{\ea}	& \text{if $l$ ramifies and $\e=\ea+\epsilon_E+1$},
	\\
	0	& \text{if $l$ ramifies and $\e>\ea+\epsilon_E+1$}.
  \end{array}
\right.
$$
\end{list}
	\qed
\end{thm}

\begin{cor}
If $E/\fp$ is ordinary, then $N_{E, p}(l^\e)$ depends only on $\ea$ and
$
\epsilon_E
$.

	\qed
\end{cor}

\begin{remark}
Ito also determines $N_{E, p}(l^\epsilon)$ when $E/\fp$ is supersingular;
again the answer depends only on $\ea$ and
$
\epsilon_E
$.
We do not reproduce the result here since the supersingular case will not
arise in our application.  Also, note that Ogg's formula (\ref{ogg}) is in
agreement with (case (iii) of) Ito's Theorem.
\end{remark}

The Corollary suggests a new notation.  Given $a$ and $p$ as above, let
$f>0$ be an integer with
$
f^2 | (a^2-4p)
$
and 
$
(a^2-4p)/f^2 \equiv\mymod{0, 1}{4}
$.
Then $\Q(\pi_a)$ has a unique order of discriminant
$
(a^2-4p)/f^2
$
which contains $\Z[\pi_a]$, and which, by Lemma \ref{lem:waterhouse},
corresponds to the endomorphism ring of some elliptic curve $E/\fp$
with $p+1-a$ points.  If $E, E'$ are two such curves, Ito's Theorem says that
$
N_{E, p}(M) = N_{E', p}(M)
$
for all $M$.  Thus
$$
N_{a, f, p}(l^\epsilon) := N_{E, p}(l^\epsilon)
$$
is well-defined.  We now come to the main result of this section.

\begin{lem}
	\label{lem:mine}
Let $M>0$ be odd.
Set
$
M_f := (M, f)
$.
Then we have the equality
$$
N_{a, f, p}(M)
=
\frac{\Psi(M)}{\Psi(M/ M_f)} \sigma(a, f, p, M),
$$
where
$
\Psi(m) := m\prod_{q|m}(1+1/q)
$,
$q$ runs through the prime divisors of $m$; and
$$
\sigma(a, f, p, M)
=
\#\{ \mymod{x}{M}: x^2 - ax + p \equiv\mymod{0}{M_f M} \}.
$$
\end{lem}

\vsp

\begin{remark}
Note that the $x$'s in $\sigma(a, f, p, M)$ are well-defined: since $M_f$
divides the conductor of $\Z[\pi_a]$, if
$
h_{a, p}(\pi)\equiv\mymod{0}{M_f M}
$
then $\pi$ is a double root of $\mymod{h}{M_f M}$, whence
$
h'_{a, p}(\pi) \equiv\mymod{0}{M}
$.
Thus for any
$
\lambda\equiv\mymod{0}{M}
$
we have
$
h_{a, p}(\pi+\lambda)
=
h_{a, p}(\pi) + h_{a, p}'(\pi)\lambda + \lambda^2
\equiv
\mymod{0}{M_f M}
$,
as desired.

\end{remark}

\begin{proof}
It suffices to assume that $M$ is an odd-prime power $l^\e$, and $\epsilon_E$
in Ito's Theorem is equal to
$$
\epsilon_f := \text{ord}_l(c_a/f).
$$
The two cases $1\le \ea - \ef$ and $\ea - \ef < \e \le \ea + \ef$ are
immediate consequence of the first part of the following Lemma.

\begin{lem}
	\label{lem:poly}
Let $g\in\Z[x]$ be monic and quadratic, and let $l>2$ be a prime.

{\rm (a)}
Suppose
$
l^{2m_0} | \text{disc}(g)
$.
Then for any integer $1\le m\le m_0$,
$$
\#\{ \mymod{x}{l^m}:  g(x)\equiv\mymod{0}{l^m} \}
=
l^{m_0 - [(m+1)/2]}.
$$
Moreover, there exists a unique $\mymod{x}{l^{m_0}}$ such that
$
g(x)\equiv\mymod{0}{l^{2m_0}}
$.

{\rm (b)}
If $l ||\text{disc}(g)$, then $g(x)$ has a unique solution
$
\mymod{}{l}
$
and has no solution $\mymod{}{l^n}$ for any $n>1$.

\end{lem}

\begin{proof}
Write $g(x) = x^2 - ax + b$.  Since $l>2$,
$$
g(a/2) = -(a^2 - 4b)/4 \equiv\mymod{0}{l^{2m_0}}
$$
by hypothesis.  Furthermore, the quadratic formula shows that $a/2$ is a double
root of
$
\mymod{g}{l^{m_0}}
$,
whence
$
g'(a/2) \equiv\mymod{0}{l^{m_0}}
$.
Thus for any $\beta\in\Z$ and any $m\le m_0$,
\begin{equation}
g(a/2+\beta) = g(a/2) + g'(a/2)\beta + \beta^2 \equiv \mymod{\beta^2}{l^m}
	\label{modm}
\end{equation}
Thus
$
g(a/2+\beta) \equiv\mymod{0}{l^m}
$
if and only if
$
\beta\equiv\mymod{0}{l^{[(m+1)/2]}}
$.
This gives the first part of (a).  To get the second part, note that if
$
\beta\equiv\mymod{0}{l^{[(m_0+1)/2]}}
$,
then combining 
$
l^{2m_0} | g(a/2)
$
and
$
l^{m_0} | g'(a/2)
$
with (\ref{modm}), we get
$
g(a/2+\beta)  \equiv \mymod{0}{l^{m_0 + [(m_0+1)/2]}}
$
if and only if 
$
\beta\equiv\mymod{0}{l^{[(3m_0+1)/4]}}
$.
Repeat this argument finitely many times and we get
$
g(a/2+\beta)  \equiv \mymod{0}{l^{2m_0}}
$
if and only if 
$
\beta\equiv\mymod{0}{l^{m_0}}
$,
as desired.

As for Part (b), $l|\text{disc}(g)$ implies that $a/2$ is also a repeated
root $\mymod{0}{l}$, whence $l| g'(a/2)$ as well, and (\ref{modm}) shows that
$
g(x)\equiv\mymod{0}{l^2}
$
has at most one solution.  But
$
g(a/2) = -(a^2-4b)/4
$
is exactly divisible by $l$, by hypothesis.

\end{proof}

Returning to the proof of the last  case
$
\e > \ea + \ef
$
 of Lemma \ref{lem:mine},
denote by $\Z[\pi_0]$ the maximal order in $\Q(\pi_a)$, and by $g_0\in\Z[x]$
the minimal polynomial of $\pi_0$.  Then
$
\pi_a = b+ \lambda \pi_0
$
for some integers $b, \lambda$ with $l^{\ea}|| \lambda$, so with the
change of variable
$
z = (x-b)/\lambda
$,
$$
g(x)\equiv\mymod{0}{l^{\e+\ea-\ef}}
	\Longleftrightarrow
g_0(z)\equiv\mymod{0}{l^{\e+\ea-\ef-2\ea}}.
$$
Thus
\begin{eqnarray*}
&&
\#\{ \mymod{x}{l^\e}: \hspace{17.5pt} g(x) \equiv\mymod{0}{l^{\e+\ea-\ef}} \}
\\
&=&
\#\{ \mymod{z}{l^{\e-\ea}}: g_0(z) \equiv\mymod{0}{l^{\e+\ea-\ef-2\ea}} \}.
\end{eqnarray*}
If
$
l\nmid \text{disc}(g_0)
$,
then the lifting argument in the proof of Hensel's Lemma shows that
the number of solutions 
$
\mymod{}{l^{\e-\ea}}
$
of
$
g_0(z) \equiv\mymod{0}{l^{\e+\ea-\ef-2\ea}}
$
is either $2$ or $0$, depending on whether $l$ is split or inert in
$
\Q(\pi_a)
$.
This gives Lemma \ref{lem:mine} when
$
\epsilon>\ea-\ef$ and $l\nmid \text{disc}(g_0)
$.
The last subcase
$
l|\text{disc}(g_0)
$
is handled by Lemma \ref{lem:poly}(b).  This completes the proof of
Lemma \ref{lem:mine}.

\end{proof}

We close this section with a computation we need for the next Section.
Fix
$
\pi, \ov{\pi}\in\C
$
so that $\pi\ov{\pi} = p\in\Z$.  For any integer $k\ge 0$, define
$$
Q(\pi,k) := \frac{\pi^{2k+1}-\ov{\pi}^{2k+1}}{\pi - \ov{\pi}}.
$$
Then
$
Q(\pi,k) = (\pi^{2k} + \ov{\pi}^{2k})
+
\text{$\Z$-linear combinations of $(\pi^{2i}+\ov{\pi}^{2i})$ for $0\le i< k$}
$,
whence
$
(\pi+\ov{\pi})^{2n}
$
is a $\Z$-linear combinations of the $Q(\pi,k)$.  Our goal is to determine
these coefficients.

For any $n\ge 1$, write
$
(\pi+\ov{\pi})^{2n}
=
\sum_{i=0}^n c_p(n, i) Q(\pi,i)
$.
Using the identities
$
(\pi+\ov{\pi})^2 = Q(\pi,1) + p
$
and, for any $k>0$,
$
(\pi+\ov{\pi})^2 Q(\pi,k) = Q(\pi,k+1) + 2p Q(\pi,k) + p^2 Q(\pi,k-1)
$,
we readily obtain the recurrence relations
\begin{eqnarray*}
c_p(n+1, n+1) & = & 1,
\\
c_p(n+1, n) & = & c_p(n, n-1) + 2p,
\\
c_p(n+1, i) & = & p^2 c_p(n, i+1) + 2p c_p(n, i) + c_p(n, i-1)
	\:\:
	\text{ for $1\le i< n$,}
\\
c_p(n+1, 0) & = & p^2 c_p(n, 1) + p c_p(n, 0).
\end{eqnarray*}
A simple induction then yields the desired formula: for $i<n$,
\begin{equation}
c_p(n, i) = (2i+1)\frac{(2n)!}{(n-i)!(n+i+1)!}p^{n-i}.
	\label{birch}
\end{equation}
Note that every $c_p(n, i)$ is positive.

\section{Elliptic modular surfaces}
	\label{sec:modular}

%
%

\begin{proof}[Proof of Theorem \ref{thm:modular}]
Denote by
$
f_M: \E\rarr X_0(M)
$
the elliptic modular surface associated to $X_0(M)$.  Both $f_M$ and $\E$
are defined over $\Q$.  Shioda
	\cite{shioda:modular}
shows that $f_M$ is a semistable fibration, and that its group of section is
finite even after we extend the base field to $\C$.  To simplify the notation,
write $Y$ and $X$ for $Y_0(M)$ and $X_0(M)$, respectively.  Then
\begin{eqnarray*}
\sum_{t\in X(\F_p)} a_p(\E_t)^{n}
&=&
\sum_{a^2<4p} a^{n} \sum_{\substack{t\in X(\F_p)\\ a_p(\E_t)=a}} 1
\\
&=&
\sum_{a^2<4p} a^{n} \sum_{\substack{t\in Y(\F_p)\\ a_p(\E_t)=a}} 1
+
O_X(1)
\end{eqnarray*}
since
$
|a_p(\E_t)|\le 1
$
for every bad fiber $\E_t$.   The number of
$
t\in Y(\F_p)
$
with $j(\E_t)=0$ or $1728$ is at most
$
\deg (X_0(M)\rarr X_0(1)) \ll M
$,
so 
\begin{eqnarray}
\sum_{t\in X(\F_p)} a_p(\E_t)^{n}
&=&
\sum_{a^2<4p} a^{n}
	\sum_{\substack{j\in\F_p\\ j\not=0, 1728}}
	\sum_{\substack{t\in Y(\F_p)\\ a_p(\E_t)=a\\j(\E_t)=j}} 1
+
O_X(p^{n/2}).
	\label{inner}
\end{eqnarray}
The supersingular fibers do not contribute to the sum, so for the rest of this
section, we will assume that
\begin{equation}
j(\E_t)\not=0, 1728	\text{ and } a_p(\E_t)\not=0.
	\label{assume}
\end{equation}
Thanks to Lemma \ref{lem:ogg}, the inner-most sum in
(\ref{inner}) is simply
$$
\frac{1}{2}N_{E_j, p}(M) = \frac{1}{2}N_{a, f, p}(M),
$$
where $E_j/\fp$ is any curve with $j(E_j)=j$, and $f=c_a/c_{E_j}$ as before.
As we run through all $j\in\fp$ with
$
j\not=0, 1728
$
and
$
a_p(E_j) = a
$,
by Lemma \ref{lem:waterhouse} we run through precisely all integers $f>0$
such that
\begin{equation}
\text{$f^2$ divides $a^2-4p$ and $(a^2-4p)/f^2 \equiv\mymod{0, 1}{4}$.}
	\label{f}
\end{equation}
The number of $E/\fp$ with $E(\fp)=p+1-a$ and with
$
\text{End}_{\fp}(E)
$
equals to the order in $\Q(\pi_a)$ with discriminant
$
(a^2-4p)/f^2
$,
is equal\footnote{this
	is
	due to Waterhouse \cite[Thm.~4.5]{waterhouse}.  However, note the
	correction in \cite[p.~194]{schoof:cubic}, which does not affect us.}
to
$
h( (a^2-4p)/f^2 )
$,
where
\begin{eqnarray*}
h(\Delta)
&:=&
\text{the class number of the imaginary quadratic order of discriminant $\Delta$.}
\end{eqnarray*}
Combine all these together, invoke Lemma \ref{lem:mine} and we can rewrite
(\ref{inner}) as
\begin{eqnarray}
&&
\frac{1}{2}
\sum_{a^2<4p} a^{n}
	\sum_f
	h\Bigl(\frac{a^2-4p}{f^2}\Bigr) N_{a, f, p}(M)
+
O_X(p^{n/2})
	\nonumber
\\
&=&
\frac{1}{2}
\sum_{a^2<4p} a^{n}
	\sum_f
	h\Bigl(\frac{a^2-4p}{f^2}\Bigr) 
	\frac{\Psi(M)}{\Psi(M/ M_f)} \sigma(a, f, p, M)
+
O_X(p^{n/2}),
	\label{almost}
\end{eqnarray}
where $f$ runs through all integers as in (\ref{f}) and $M_f = (M, f)$ as
in the last Section.   Since
$
\sigma(a, f, p, M) = \sigma(-a, f, p, M)
$,
if $n$ is odd then the double sum in (\ref{almost}) is zero, whence the
entire expression (\ref{almost}) is
$
O_X(p^{n/2})
$.
Substitute this back into (\ref{inner}) and we get, for odd $n$,
\begin{eqnarray*}
\sum_p \frac{-\log p}{p^{s+[n/2]+1}}  \sum_{t\in X(\F_p)} a_p(\E_t)^{n}
&\ll&
\sum_p \frac{-\log p}{p^{s+[n/2]+1}} p^{n/2}.
\end{eqnarray*}
Theorem \ref{thm:modular}(a) then follows.  To handle even $n$, we rewrite
(\ref{almost}) in terms of the Selberg trace formula for Hecke operators on
cusp forms for
$
\Gamma_0(M)
$;
Theorem \ref{thm:modular} then follows from the Weil conjecture estimate for
these traces.

\if 3\
{
We begin with the case $n>2$.  Ito's Theorem \ref{thm:ito} gives
$
N_{E, p}(M) \le 2M
$
for all $p$ and
$
\Psi(M) / \Psi(M/ M_f) \le 2M
$
as well.  If $\Delta$ is the fundamental discriminant of an imaginary
quadratic order, then
\begin{eqnarray}
h(\Delta \delta^2)
&=&
h(\Delta) \prod_{q|\Delta} \Bigl( 1 + \Bigl(\frac{\Delta}{q}\Bigr) \Bigr)
	\hspace{20pt}
	\text{by \cite[p.~????????]{cox}}
\\
&\ll_\epsilon&
|\Delta|^{1/2+\epsilon},
	\hspace{20pt}
	\text{\framebox{what if $(\delta, \Delta)>1$??}}
\end{eqnarray}
where $q$ runs through the distinct prime divisors of $\Delta$.  Thus
$$
h\Bigl(\frac{a^2-4p}{f^2}\Bigr) \frac{\Psi(M)}{\Psi(M/ M_f)}
	\sigma(a, f, p, M)
\ll_{X,\epsilon}
	(a^2-4p)^{1/2+\epsilon}
\ll_{X,\epsilon}
	a^{1+\e}.
$$
The number of divisors of an integer $n$ is
$
\ll_\epsilon n^\e
$,
so putting all these together and we see that (\ref{almost}) is bounded from
the above by
$$
\ll_{X,\e}
\:\:
\sum_{|a|<2\sqrt{p}} a^{n}  \sum_{f^2| (a^2-4p)}  a^{1+\epsilon}
\:\:
\ll_{X,\e}
\:\:
\sum_{|a|<2\sqrt{p}} a^{n+1+\epsilon}
\:\:
\ll_{X,\e}
\:\:
p^{n/2+1+\e}.
$$
Since $n>2$, this is 
$
\ll_{X,\e} p^{n-1/2+\e}
$,
whence
\begin{eqnarray*}
\raisebox{-7pt}{$\stackrel{\displaystyle\text{res}}{\scriptstyle s=1}$}\,
	\sum_p
	\raisebox{20pt}{\mbox{}}
	\frac{-\log p}
		{p^{s+n}}\!\!\sum_{t\in X_0(M)(\F_p)}
	a_p(\mathcal E_t)^{n}
&
=
&
\raisebox{-7pt}{$\stackrel{\displaystyle\text{res}}{\scriptstyle s=1}$}\,
	\sum_p
	\raisebox{20pt}{\mbox{}} \frac{-\log p}{p^{s+n}} p^{n-1/2+\e}
\:\:
=
0.
\end{eqnarray*}
}
\fi

For any fixed $p$, the number of pairs
$(a, f)$ with
$
a^2 - 4p = -3f^2
$
or
$-4f^2$ is at most four.  Since
$
N_{a, f, p}(M)
$
is bounded from the above in terms of $M$, with
$$
h_w(\Delta)
:=
\left\{
  \begin{array}{ll}
	h(\Delta)	& \text{if $\Delta<-4$}
	\\
	1/2		& \text{if $\Delta=-4$}
	\\
	1/3		& \text{if $\Delta=-3$}
  \end{array}
\right.
$$
and recall that $a=\pi_a + \ov{\pi}_a$, for even $n$ we can further rewrite
(\ref{almost}) as
\begin{eqnarray}
&&
\frac{1}{2}
\sum_{a^2<4p} (\pi_a+\ov{\pi}_a)^{n}
	\sum_f
	h_w\Bigl(\frac{a^2-4p}{f^2}\Bigr) 
	\frac{\Psi(M)}{\Psi(M/ M_f)} \sigma(a, f, p, M)
+
O_X(p^{n/2})
	\nonumber
\\
&=&
\sum_{i=0}^{n/2}
c_p(n/2, i)
\frac{1}{2}
\sum_{a^2<4p} Q(\pi_a, i)
	\sum_f
	h_w\Bigl(\frac{a^2-4p}{f^2}\Bigr) 
	\frac{\Psi(M)}{\Psi(M/ M_f)} \sigma(a, f, p, M)
+
O_X(p^{n/2}),
	\label{two}
\end{eqnarray}
where the
$
c_p(n/2, i)
$
and
$
Q(\pi_a, i)
$
are as in the end of the last Section.  The Selberg trace formula says that
the trace of the Hecke operator $T_p$ on the space of weight $2k$ modular
forms on $\Gamma_0(M)$ with
$
p\nmid M
$,
is given by  \cite[Thm.~2.2]{schoof} 
\begin{eqnarray*}
\text{trace}(T_p, S_{2k}(\Gamma_0(M)))
&=&
\frac{-1}{2}
\sum_{a^2<4p}
	Q(\pi_a, k-1)
\sum_f
h_w\Bigl(\frac{a^2-4p}{f^2}\Bigr)
	\frac{\Psi(M)}{\Psi(M/M_f)} \sigma(a, f, p, M)
\\
&&
-
\sum_{\substack{c|M\\(c, M/c)|(M,p-1)}} \phi((c, M/c))
+
\Bigl\{
  \begin{array}{ll}
	p+1 & \text{ if $2k=2$}
	\\
	0   & \text{ otherwise.}
  \end{array}
\Bigr.
\end{eqnarray*}
This allows us to rewrite (\ref{two}) as
\begin{eqnarray*}
&&
c_p(n/2, 0)p
-
\sum_{i=1}^{n/2}
	c_p(n/2, i)
	\Bigl[ \text{trace}(T_p, S_{2i}(\Gamma_0(M))) + O_X(1) \Bigr]
+
O_X(p^{n/2})
\\
&=&
\frac{p^{n/2+1} n!}{(n/2)! (n/2+1)!} 
+
O_X(p^{n/2}),
\end{eqnarray*}
since the Weil conjecture gives
$
\text{trace}(T_p, S_{2i}(\Gamma_0(M))) \ll p^{i-1/2}
$
while
(\ref{birch}) gives
$
c_p(k, i)\ll p^{k-i}
$.
Substitute these back into  (\ref{inner}) and we get Theorem \ref{thm:modular}
for even $n$.

\end{proof}

\begin{remark}
Our use of the trace formula above is inspired by the work of
Birch \cite{birch} on the asymptotic behavior of \textit{even} moments of
values of $a_p(E)$ over all elliptic curves over
$\fp$.  This is essentially the $M=1$ case above; of course,
when $M=1$ the problem of counting cyclic $M$-isogenies does not arise.

\end{remark}

\section{Tate conjecture for $\E$ and $\etwotilde$}
	\label{sec:tate}

To simplify our notation, write $\E$ for $\E_M$ and write $C$ for $X_0(M)$.
For the rest of this section the ground field is $\Q$.

We begin with $\E$.  Following the notation in the Introduction, denote by
$
T_\E
$
the subgroup of $NS(\E/\ov{\Q})$ generated by the zero section of
$
\efc
$
together with all the $\ov{\Q}$-components of the fibers of $f$.  From the
theory of modular curves we know that these irreducible components are all
defined over (subfields of) the cyclotomic field
$
\Q(\zeta_M)
$.
Consequently, the Artin $L$-function 
$
L(T_\E\otimes\Q, s)
$
attached to the $\Q[G_\Q]$-module $T_\E\otimes\Q$ is a product of Dirichlet
$L$-series.  In particular, 
$
L(T_\E\otimes\Q, s)
$
has a meromorphic continuation to $\C$ and is holomorphic for $Re(s)>0$ except
for a pole at $s=1$ of order
$
\text{rank } T_\E(\Q)
$.
On the other hand, Rosen and Silverman \cite[p.~52-53]{mike_joe} shows that,
for Re$(s)>3/4$,
$$
\sum_p \frac{-\log p}{p^{s+1}}\sum_{t\in X_0(M)(\F_p)} a_p(\E_t)
=
-\frac{d}{ds} \log L_2(\E/\Q, s+1)
+
\frac{d}{ds} \log L(T_\E/\Q, s)
+
O(1).
$$
By Theorem \ref{thm:modular}, the left side has an analytic continuation to
Re$(s)>1$.  It follows that 
$
L_2(\E/\Q, s)
$
has an analytic continuation to $Re(s)>7/4$ except for a pole at $s=2$ of
order equal to 
$
\text{rank } T_\E(\Q)
$.
Finally,  Shioda \cite{shioda:modular} shows that the group of sections of
$
\E\rarr C
$
is finite (even when the base field is $\C$), so by the Shioda-Tate theorem,
the order of pole of
$
L_2(\E/\Q, s)
$
at $s=2$ is in fact
$
\text{rank } NS(\E/\Q)
$.
This completes the proof of the Tate conjecture for $\E/\Q$.

\if 3\
{
\begin{lem}[\protect{Cox \cite[p.~320]{cox}}]
	\label{lem:cox}
Let $E, E'$  be two elliptic curves over $\F_p$ with the same $j$-invariant.
Then

{\rm (a)}
If $j(E)\not=0, 1728$, then $\#E(\fp) = \pm \#E'(\fp)$.

{\rm (b)}
The two curves are $\F_p$-isomorphic if and only if
$\#E(\fp) = \#E'(\fp)$.

	\qed
\end{lem}

\noindent
By Lemma \ref{lem:cox}, 
}
\fi

\vsp

\begin{lem}
	\label{lem:ord}
$
-\raisebox{-7pt}{$\stackrel{\displaystyle{\rm ord}}{\scriptstyle s=2}$}\,
L_2(\etwotilde/\Q, s)
\ge
{\rm rank} \bigl[ H^2_\et(\etwotilde\otimes\ov{\Q}, \Q_l(1))^{G_\Q} \bigr].
$
\end{lem}

\begin{proof}
Part of the following calculation is inspired by
	\cite[Lem.~1.6]{schoen}.

From the spectral sequence for
$
\tilde{\pi}: \etwotilde\rarr X_0(M)
$
we get a $G_\Q$-equivariant filtration on
$
H^2_\et(\etwotilde\otimes\ov{\Q}, \Q_l) =: L^0
$
with
$$
L^2 \simeq H^2(R^0),
\:\:\:
L^1/L^2 \simeq H^1(R^1),
\:\:\:
L^0/L^1 \simeq H^0(R^2),
$$
where 
$
H^i(R^j)
:=
H^i_\et(X_0(M)\otimes\ov{\Q}, R^j\tilde{\pi}_\ast \Q_{l, \etwotilde})
$.
Thus
$
L_2(\etwotilde/\Q, s)
$
is the product of the $L$-functions associated to 
$
H^{2-t}(R^t)
$
for $t=0, 1, 2$.  Moreover,
\begin{equation}
\sum_{t=0}^2  \text{rank} \bigl[ H^{2-t}(R^t)(1)^{G_\Q} \bigr]
\ge
\text{rank} \bigl[ H^2_\et(\etwotilde\otimes\ov{\Q}, \Q_l(1))^{G_\Q} \bigr].
	\label{ineq}
\end{equation}
To prove the Lemma it then suffices to show that for every $t$,
\begin{equation}
-
\raisebox{-7pt}{$\stackrel{\displaystyle\text{ord}}{\scriptstyle s=2}$}\,
L( H^{2-t}(R^t), s)
\ge
\text{rank} \bigl[ H^{2-t}(R^t)(1)^{G_\Q} \bigr].
	\label{claim}
\end{equation}

We begin with $t=0$.  Since
$
R^0 \tilde{\pi}_\ast \Q_{l, \etwotilde}\simeq \Q_{l, X_0(M)}
$,
we have
$
H^2(R^0) \simeq \Q_l(-1)
$.
In particular, 
$
\text{trace}(\text{Frob}_p, H^2(R^0)) = p
$.
Thus
$
L( H^2(R^0), s) = \zeta(s-1)
$
has a simple pole at $s=2$ while
$
\text{rank} H^2(R^0)(1)^{G_\Q} = 1
$,
so both sides of (\ref{claim}) are $1$ when $t=0$.

Denote by $Z$ the singular locus of $\efc$; it is also the singular locus
for
$
\etwotilde\rarr C
$.
Denote by $i: Z\rarr C$ the closed immersion, and by $j: U\rarr C$ the
open embedding of $U=C-Z$ in $C$.
For any $l$-adic sheaf $\mathcal F$ on $C$ we have an exact sequence
$$
0\rarr i_\ast i^! \mathcal F \rarr \mathcal F \rarr j_\ast j^\ast \mathcal F.
$$
This becomes a short exact sequence for
$
\mathcal F = R^t \tilde{\pi}_\ast \Q_{l, \etwotilde}
$,
by the local invariant cycle theorem.  Moreover,
\begin{eqnarray*}
j^\ast R^t \tilde{\pi}_\ast \Q_{l, \etwotilde}
&\simeq&
R^t {\pi_U}_\ast (j_\E^\ast)\Q_{l, \etwotilde}
				\hspace{20pt}\text{proper base change}
\\
&\simeq&
R^t {\pi_U}_\ast\Q_{l, \E^{(2)}_U},
\end{eqnarray*}
where
$
f_U: \E_U\rarr U
$
is the pull-back to $U$ of
$
f: \E\rarr C$ via $j
$
and
$
\pi_U: \E^{(2)}_U\rarr U
$
the corresponding map.  Consequently,
\begin{eqnarray}
\text{trace}(\text{Frob}_p, H^{2-t}(R^t) )
&=&
\text{trace}(\text{Frob}_p, H^{2-t}_\et(X_0(M)\otimes\ov{\Q}, i_\ast i^! R^t \tilde{\pi}_\ast \Q_{l, \etwotilde}))
+
	\label{a}
\\
&&
\text{trace}(\text{Frob}_p, H^{2-t}_\et(X_0(M)\otimes\ov{\Q}, j_\ast R^t {\pi_U}_\ast \Q_{l, \E^{(2)}_U})).
	\label{b}
\end{eqnarray}
Since $i_\ast i^!\mathcal F$ is a skyscraper sheaf, for $t=1$ the Kunneth
formula gives
\begin{eqnarray*}
\text{trace}(\text{Frob}_p, H^1(R^1) )
&=&
\text{trace}(\text{Frob}_p, H^1_\et(X_0(M)\otimes\ov{\Q}, j_\ast R^1 {\pi_U}_\ast \Q_{l, \E^{(2)}_U}))
\\
&&
\text{trace}(\text{Frob}_p, H^1_\et(X_0(M)\otimes\ov{\Q}, j_\ast R^1 {f_U}_\ast \Q_{l, \E_U}))^{\oplus 2}.
\end{eqnarray*}
Deligne \cite{deligne} shows that
$
H^1_\et(X_0(M)\otimes\ov{\Q}, j_\ast R^1 {f_U}_\ast \Q_{l, \E_U})
$
is canonically identified with the space of weight $3$ modular forms on
$\Gamma_0(M)$ plus its complex conjugate.  There are no weight $3$ modular
forms on $\Gamma_0(N)$,  so for $t=1$ both sides of (\ref{claim}) are zero.

Finally, consider the case $t=2$.  The Kunneth formula gives
\begin{eqnarray*}
\lefteqn{H^0_\et(X_0(M)\otimes\ov{\Q}, j_\ast R^2 {\tilde{\pi}}_{U\ast} \Q_{l, \E^{(2)}_U})}
\\
&\simeq&
H^0_\et(X_0(M)\otimes\ov{\Q}, j_\ast
	\underbrace{R^2 \tilde{f_U}_\ast\Q_{l, \E_U}}_{\Q_l(-1)}
)^{\oplus 2}
\oplus
H^0_\et(X_0(M)\otimes\ov{\Q}, j_\ast
	\underbrace{(R^1 \tilde{f_U}_\ast\Q_{l, \E_U})^{\otimes 2}}_{\Q_l(-1)}
).
\\
&\simeq&
\Q_l(-1)^{\oplus 2}\oplus \Q_l(-1),
\end{eqnarray*}
so (\ref{b}) contributes $3$ to both sides of (\ref{claim}).  As for (\ref{a}),
$$
H^0_\et(X_0(M)\otimes\ov{\Q}, i_\ast i^! R^2 \tilde{\pi}_\ast \Q_{l, \etwotilde})
\simeq
\oplus_x H^0_\et(x\otimes\ov{\Q}, {i_x}_\ast i^! R^2 \tilde{\pi}_\ast \Q_{l, \etwotilde}),
$$
where $x$ runs through the closed points of $Z\subset C$ and 
$
i_x: x\rarr C
$
denotes the closed immersion.  From the theory of modular curves we see that
$$
k(x) := \text{residue field of $x$}
$$
is a subfield of the cyclotomic field $\Q(\zeta_M)$.  The fibers of $f$ are
also defined over some subfields of $\Q(\zeta_M)$, and hence the same holds
for the fibers of $\tilde{\pi}$.  That means as $G_{k(x)}$-modules,
$
{i_x}_\ast i^! R^2 \tilde{\pi}_\ast \Q_{l, \etwotilde}
$
is a direct sum of Abelian characters $\chi$'s.  The $L$-function associated to
$
H^0_\et(x\otimes\ov{\Q}, {i_x}_\ast i^! R^2 \tilde{\pi}_\ast \Q_{l, \etwotilde})(1)
$
is therefore a product of Abelian $L$-series
$
L(s, \chi)
$.
It is classical
	that
	at\footnote{we
	shift from $s=2$ to $s=1$ because of the Tate twist in
	$
	H^0_\et(x\otimes\ov{\Q}, {i_x}_\ast i^! R^2 \tilde{\pi}_\ast \Q_{l, \etwotilde})(1)
	$}
$s=1$,  $L(s, \chi)$ either has a pole or has neither a pole nor a zero;
we pick up a pole precisely when
$
H^0_\et(x\otimes\ov{\Q}, \chi)(1)
$
is a $1$-dimensional $G_\Q$-fixed space of
$
H^0_\et(x\otimes\ov{\Q}, {i_x}_\ast i^! R^2 \tilde{\pi}_\ast \Q_{l, \etwotilde})(1)
$.
Thus the contribution from
(\ref{b})
to the two sides of (\ref{claim}) are equal as well, and the Lemma
follows.

\end{proof}

\if 3\
{
\begin{equation}
\text{trace}(\text{Frob}_p, H^2_\et(\etwotilde\otimes\ov{\Q}, \Q_l) )
=
\sum_{t=0}^2
\text{trace}(\text{Frob}_p, H^{2-t}_\et(X_0(M)\otimes\ov{\Q}, R^t \tilde{\pi}_\ast \Q_{l, \etwotilde})).
	\label{trace}
\end{equation}
Thus
$
L_2(\etwotilde/\Q, s)
$
is the product of the $L$-function associated to each
$
H^{2-t}_\et(X_0(M)\otimes\ov{\Q}, R^t \tilde{\pi}_\ast \Q_{l, \etwotilde})
$.
We now examine each of the terms on the right side.
}
\fi

\begin{remark}
The Lemma also holds for the universal elliptic curve over
$
X_1(M)/\Q(\zeta_M)
$.
This time
$
H^1(R^1)
$
is not empty; but (\ref{claim}) remains true since both sides are still zero.
For the analytic side it follows from a deep non-vanishing theorem
	\cite{jacquet}.
For the cohomological side, this follows from the fact that no power of the
cyclotomic character is a quotient of the Galois representation associated to
a weight $3$ cusp forms on $\Gamma_1(M)$.
\end{remark}

\vsp

We now return to Tate's conjecture for $\etwotilde$.  As in
	(\ref{firstterm}) --- (\ref{thirdterm}),
we get
\begin{eqnarray}
\lefteqn{\frac{d}{ds}\log L_2(\etwotilde/\Q, s)}	\nonumber
\\
&=&
2 \frac{d}{ds}\log L_2(\E/\Q, s)
-
\sum_p
\frac{-\log p}{p^{s}}(1+b_2)
-
\sum_p
\frac{-\log p}{p^{s+2}}
	\sum_{t\in C(\F_p)} a_p(\E_t)^2.
\label{tatefortwo}
\end{eqnarray}
Invoke Theorem \ref{thm:modular} and the analytic Tate conjecture just proved
for $\E/\Q$, we see that
$
L_2(\etwotilde/\Q, s)
$
has an analytic continuation to $Re(s)>7/4$ except possibly for a pole at
$s=2$.
Moreover, compute the residues on both sides of (\ref{tatefortwo}) at $s=1$
and we get
\renewcommand{\arraystretch}{1.2}
$$
\begin{array}{lllll}
\text{rank} NS(\etwotilde/\Q)
&\ge&
2 \text{ rank} NS(\E/\Q) + b_2
&&
\text{by Theorem \ref{thm:divisor}}
\\
&=&
-\raisebox{-7pt}{$\stackrel{\displaystyle\text{ord}}{\scriptstyle s=2}$}\,
L_2(\etwotilde/\Q, s)
&&
\text{by (\ref{tatefortwo})}
\\
&\ge&
\text{rank} H^2_\et(\E\otimes\ov{\Q}, \Q_l)^{G_\Q}
&&
\text{by Lemma \ref{lem:ord}}
\\
&\ge &
\text{rank} NS(\etwotilde/\Q)
\end{array}
\renewcommand{\arraystretch}{1}
$$
since
$
NS(\etwotilde/\Q)\otimes\Q_l
			\hookrightarrow
H^2_\et(\E\otimes\ov{\Q}, \Q_l)^{G_\Q}
$.
Thus we have equalities all the way, whence the order of pole at $s=2$ of
$
L_2(\etwotilde/\Q, s)
$
is 
$
\text{rank} NS(\etwotilde/\Q)
$,
as desired.

	\qed

\begin{remark}
	\label{remark:geom}
From this analytic proof of Tate's conjecture for $\etwotilde/\Q$
we see that the inequality in Lemma \ref{lem:ord} is in fact an equality.
That means
\begin{itemize}
\item
(\ref{ineq}) is in fact an equality, and
\item
in the notation of the proof of the Lemma, for each
$
x\in Z
$,
the $G_\Q$-fixed subspace of
$
H^0_\et(x\otimes\ov{\Q}, {i_x}_\ast i^! R^2 \tilde{\pi}_\ast \Q_{l, \etwotilde})(1)
$
is spanned by precisely to the divisors $A_{ij}(t)$ in (\ref{list1}) for the
cusp $x$ plus the blowup divisors.
\end{itemize}
With a more delicate analysis we should be able to prove these two statements
directly.  Furthermore, the contribution from the term
$
H^2(R^0)(1)^{G_\Q}
$
corresponds to a vertical fiber of
$
\widetilde{\E_M}^{(2)}
$,
while the contributions from
$
H^0_\et(X_0(M)\otimes\ov{\Q}, j_\ast R^2 {\tilde{\pi}}_{U\ast} \Q_{l, \E^{(2)}_U})
$
correspond to the divisors
$
\tilde{E}_1, \tilde{E}_2
$
and
$
\tilde{\Delta}
$
in (\ref{list2}).  Putting these together and we should get a purely geometric
proof --- modulo analytic and Galois properties of weight $3$ modular forms --- of Tate's conjecture
for
$
\widetilde{\E_M}^{(2)}
$.
We present the analytic argument here because of its own interest, and
because it serves as a guide on how to prove Tate's conjecture
for $\etwotilde$ for a general semistable fibration
$
\E\rarr C
$.
Specifically, the analytic proof above depends on three ingredients:
\begin{itemize}
\item
Tate's conjecture for $\E_M$,
\item
the crucial Lemma \ref{lem:ord}, and
\item
the Galois module structure of the bad fibers.
\end{itemize}
To extend our proof of Lemma \ref{lem:ord} to a general semistable fibration
$
\efc
$,
we need to understand the monodromy representation
$
H^1(R^1)
$.
This is also the main obstacle to proving Tate's conjecture for $\E$.  For the
universal elliptic curve over the Fermat curve
(cf.~Remark \ref{remark:fermat}), we expect the complex multiplication
structure of the Fermat curve to greatly facilitate the monodromy computation.
In additional, the Galois module structure of the bad fibers of 
Fermat fibrations is well understood.  These issues are currently under
investigation.

\end{remark}

\vsp

\section{Isotrivial fibrations}
	\label{sec:isotrivial}

\if 3\
{
To gain some insight into this question, we used the computer algebra package
PARI to evaluate this residue for five surfaces: the elliptic modular
surface over the modular curves $X_1(5), X_1(7)$ and $X_1(11)$;
as well as the surfaces defined by
$
y^2 = x(x-1)(x-t)
$
and by
$
y^2 = x(x-1)(x-t(2-t))
$.
The last one has rank $1$ over $\Q(t)$, while the other four all have rank
$0$ over the function fields of the corresponding base curves.  Using the first
$100,000$ primes, in all five cases the residues converge to $1$
to four decimal places, albeit slowly (not too surprising perhaps, since as
we will see in the next subsection, we are essentially dealing
with the logarithmic derivative of a function at a point where the function
has a pole of finite order).

The five examples above are non-isotrivial, i.e.~the $j$-invariants
of the smooth fibers are not constant.  For isotrivial surfaces we can
evaluate this sum in many cases.
}
\fi

\begin{proof}[Proof of Theorem \ref{thm:isotrivial}]

Let $f: \E\rarr C$ be a non-split, isotrivial elliptic fibration over $\Q$.
Then we can find an integer
$
n\in\{2, 3, 4, 6\}
$,
a smooth curve $E/\Q$ of genus
	one\footnote{$f$
	is not assumed to have a section, so $E/\Q$ need not be an elliptic
	curve.}
and a non-constant function
$
g\in \Q(C)
$
such that, for all but finite many points
$
t\in C(\ov{\Q})
$,
the fiber $f^{-1}(t)$ is the $n$-th order twist of $E/\Q$ by $g(t)$; in other
words,
$
f^{-1}(t)
$
is isomorphic to $E$ over the extension
$
\Q(t, \sqrt[n]{g(t)})
$.

We begin with $n=2$, i.e.~a quadratic twist family.  Then
$
a_p(\E_t) = \pm a_p(E)
$
for every smooth fiber $\E_t$, whence by the Weil conjecture,
\begin{eqnarray}
\sum_p
    \frac{\log p}{p^{s+2}}   \sum_{t\in C(\F_p)} a_p(\E_t)^2
&=&
\sum_p    \frac{\log p}{p^{s+1}} a_p(E)^2
+ 
O_\E\Bigl(\sum_p \frac{\log p}{p^{s+3/2}} a_p(E)^2\Bigr).
	\label{rankin-selberg}
\end{eqnarray}
Note that up to a term that is holomorphic for
$
Re(s)>1/2
$,
the first term on the right side is the logarithmic derivative of the
Rankin-Selberg convolution of the $L$-series of the Jacobian $J_E/\Q$ of
$E/\Q$.  Since $J_E/\Q$ is a \textit{modular} elliptic curve, the
convolution $L$-series has a simple pole at $s=0$ and $1$ and is holomorphic
for
$Re(s)>1$ \cite{li}, whence the residue of (\ref{rankin-selberg}) at
$s=1$ is exactly one.  This gives Theorem \ref{thm:isotrivial}(a) for 
quadratic twist families.

Next, consider the case $n=4$.  Then $E/\Q$ is a (possibly trivial)
$\Q$-torsor of 
$
E_0: y^2 = x^3 + x
$,
and $g\in k(C)$ is not a $4$-th power.  Every smooth curve of genus one over
any finite field $\F$ has a $\F$-rational point, so
$
E/\F_p
$
is $\F_p$-isomorphic to
$
E_D: y^2 = x^3+D_p x
$
for some element $D_p\in\F_p$.  Consequently, except for a finite
number of
$
t\in C(\F_p)
$
(the number of which is bounded from the above independent of $p$),
$
\E_t/p
$
is isomorphic to
$
E_{D_p g(t)}/\F_p
$.
It is classical \cite{ireland-rosen}  that for $p$ sufficiently large,
$$
a_p(\E_{D_p g(t)})
=
\Bigl\{
  \begin{array}{lll}
  0	&&	\text{if $p\equiv\mymod{3}{4}$}
  \\
  \chi_p( D_p g(t))^2 a_p(E)	&&	\text{if $p\equiv\mymod{1}{4}$,}
  \end{array}
\Bigr.
$$
where for
$
p\equiv\mymod{1}{4}
$,
$
\chi_p:(\Z/p)^\times\rarr\C
$
is a character of order $4$.  Consequently,
\begin{eqnarray}
\sum_p
    \frac{\log p}{p^{s+2}}   \sum_{t\in C(\F_p)} a_p(\E_t)^2
&=&
\Bigl(
   \sum_{p\equiv 3(4)} \frac{\log p}{p^{s+2}} a_p(E)^2
\Bigr)
	\cdot O_\E(1)
+
	\label{04}
\\
&&
\sum_{p\equiv 1(4)}
    \frac{\log p}{p^{s+2}} a_p(E)^2 \chi_p(D_p)^2 \sum_{t\in C(\F_p)}
						\chi_p(g(t))^2.
	\label{14}
\end{eqnarray}
Denote by
$
C'/\F_p
$
the double cover of $C$ given by the (possibly singular) model
$
z^2 = g(t)
$.
Then the inner sum in (\ref{14}) is simply
$$
\sum_{t'\in C'(\F_p)} a_p(C') + O_\E(1),
$$
where the $O$-constant is independent of $p$.  This is
$
O_\E(\sqrt{p})
$
by the Weil conjecture, so the residue at $s=1$ of the left side of (\ref{04})
is zero.  The same
argument applies \textit{mutatis mutandis} to the case $n=3$ and $n=6$.

\end{proof}

\begin{remark}
In order to apply this argument to general isotrivial surfaces over a number
field, we need to understand the analytic properties for the convolution of
$L$-series
of elliptic curves over number fields.  That seems to be a very difficult
problem.  

Using known results on symmetric third- and fourth-powers $L$-functions,
we can extend the argument above to $a_p^3$- and $a_p^4$-sums of isotrivial
fibrations over $\Q$.  Higher $a_p$-powers are beyond current techniques.
\end{remark}

\if 3\
{
These considerations together lead us to the following 

\begin{quest}
Let $\efc$ be a non-split elliptic surface.  Is it true that
$$
\underset{s=1}{\text{ res }}
  \sum_p
    \frac{\log p}{p^{s+2}}   \sum_{t\in C(\F_p)} a_p(\E_t)^2
\stackrel{?}{=}
\biggl\{
  \begin{array}{lll}
  0	&	\text{if $\efc$ is a cubic, quartic }
  \\
	&	\text{or sextic twist family}
  \\
  1	&	\text{otherwise}
  \end{array}
\biggr.
?
$$
\end{quest}

It would be of great interest to extend the argument above to
cover all cubic, quartic and sextic twist families over a general base curve
$C$ over $\Q$.  This requires working with character sum estimates over
a smooth projective curve.  Results of Katz \cite{katz} should be useful
here.

\begin{remark}
In this  subsection we approach the $a_p(\E_t)^2$-sum from an analytic
point of view.  In the next subsection we revisit this sum, and especially
the question as to whether or not it is $1$, from the geometric standpoint.
The geometric construction there also suggests that 
there is one important class of non-isotrivial elliptic fibrations for which
we could settle Question 3, i.e.~the class of elliptic
modular surfaces over modular curves.
\end{remark}

}
\fi

\if 3\
{

\begin{remark}
	\label{remark:four}
We claim that the residue (????) cannot be four; otherwise, the set of $\p$
with $|a_\p(\E_t)|< 2N_\p$ must have density zero; that is impossible.

\framebox{WHY??????}
\end{remark}
}
\fi

\if 3\
{
\section{Symmetric-square fibrations}
	\label{sec:symmetric_square}

Every fiber-square construction has an accompanying symmetric-square
quotient.  We now investigate that.
First, let
$
\E_1\stackrel{f_1}{\rarr} C
$
and
$
\E_2\stackrel{f_2}{\rarr} C
$
be two semistable fibrations over the same curve $C$; the singular locus in
$C$ of the $f_i$ need not be common.  Suppose
$
S\stackrel{\pi}{\rarr}C
$
is a genus two fibration whose associated Jacobian fibration is isogenous
to
$$
f_1\times_C f_2: \E_1\times_C\E_2 \rarr C.
$$
Repeat the calculation in Section \ref{sec:power} for this genus
two fibration and we get
\begin{eqnarray*}
\sum_{t\in C(\F_\p)} a_\p(S_t)^2
&=&
(N_\p^2+N_\p)
\times
\bigr(
  1+ (\#\text{singular points of $\stwo$})
\bigr)
\\
&&
+
	N_\p\bigl(
	  t_{\p, 2}( \stwotilde )  - 2b_\p(S_t)
	\bigr)
+
	O_{S, C}(N_\p^{3/2}).
\end{eqnarray*}
The argument in \cite[Lem.~1]{schoen} readily shows that
$
a_\p(\stwotilde) = a_\p(C)
$.
Apply Tate's conjecture as before and we get
\begin{eqnarray}
\text{rank }NS( \stwotilde/k )
&
\stackrel{\text{Tate}_?}{=}
&
1
+
(\#\text{singular points of $\stwo$})
+
2 \text{ rank $(\text{Jac }S)(K_C)$}
	\label{ns_s2}
\\
&&
\hspace{8pt}
+
\underset{s=1}{\text{ res }}
  \sum_\p
    \frac{\log N_\p}{N_\p^{s+2}}   \sum_{t\in C(\F_\p)} a_\p(S_t)^2.
\nonumber
\end{eqnarray}
It follows from 
\cite[p.~80]{cassels} that that for any smooth fiber
$S_t$,
\begin{equation}
a_{\p}(S_t)^2
=
a_{\p^2}(S_t) + 4N_\p - 2(N_\p+1)(a_\p(\E_{1,t}) + a_\p(\E_{2,t}) - a_\p(S_t))
	+ 2a_\p(\E_1) a_\p(\E_2).
	\label{ap2}
\end{equation}

\framebox{what about the bad fibers?!!!!!}

By the Weil conjecture, all but the last term on the right side of 
(\ref{ap2}) is
$
\ll_{S, C} N_\p^{3/2}
$,
so
\begin{eqnarray}
\underset{s=1}{\text{ res }}
  \sum_\p
    \frac{\log N_\p}{N_\p^{s+2}}   \sum_{t\in C(\F_\p)} a_\p(S_t)^2
&=&
2\underset{s=1}{\text{ res }}
  \sum_\p
    \frac{\log N_\p}{N_\p^{s+2}}   \sum_{t\in C(\F_\p)} a_\p(\E_{1, t})a_\p(\E_{2, t}).
	\label{aps2}
\end{eqnarray}
Now, take the two elliptic fibrations $\E_i\rarr C$ to be
the same. Then with $\E :=\E_i$, we have
$
\text{rank $(\text{Jac }S)(K_C)$} = 2\text{ rank }\E(K_C).
$
Combine Theorem \ref{thm:e2}, (\ref{ns_s2}) and (\ref{aps2}),
we get 
\begin{eqnarray}
\text{rank }NS( \etwotilde/k )
&
\stackrel{\text{Tate}_?}{=}
&
{\textstyle\frac{1}{2}}
\text{rank }NS( \stwotilde/k )
+
(\#\text{singular points of $\etwo$})
	\label{onemore}
\\
&&
+\:
{\textstyle\frac{1}{2}}
-
{\textstyle\frac{1}{2}} (\#\text{singular points of $\stwo$}).
	\nonumber
\end{eqnarray}

}
\fi

\section*{Acknowledgement}

I am grateful to the support and encouragement of Barry Mazur, Michael Rosen
and Joe Silverman, over the years and for this project in particular.
David Cox offers excellent tutorials on toric geometry.
%
%
%
%
Dan Abramovich patiently explains to me the fine points of Cartier
divisors.
Peter Sarnak draws our attention to the work of Birch and suggests
that it could be related to the $a_p^n$-Nagao sums.  Brian Conrad and
Eyal Markman patiently answer many of my questions.
My thanks to all of them.
Any error that remains is due to me.


\thispagestyle{empty}

\bibliographystyle{amsalpha}

\end{document}